\documentclass[12pt]{article}

\usepackage{amsmath,amsthm,amsfonts,amssymb}
\usepackage{graphicx}
\usepackage[usenames]{color}
\usepackage{hyperref}

\newtheorem{theorem}{Theorem}[section]
\newtheorem{maintheorem}{Theorem}[]
\newtheorem{lemma}[theorem]{Lemma}
\theoremstyle{definition}

\newtheorem{corollary}[theorem]{Corollary}

\newtheorem{definition}[theorem]{Definition}

\newtheorem{remark}[theorem]{Remark}
\newtheorem{proposition}[theorem]{Proposition}

\newtheorem{procedure}{Procedure}[section]

\newenvironment{evf}{\noindent\color{green} \bf EF: }{}

\renewcommand{\proof}{{\bf Proof. }}
\newcommand{\eproof}{$\bullet$}

\newcommand{\fracd}[2]{{\displaystyle\frac{#1}{#2}}}
\newcommand{\abs}[1]{\vert #1 \vert}

\newcommand{\bi}{\begin{itemize}}
\newcommand{\ei}{\end{itemize}}
\newcommand{\ben}{\begin{enumerate}}
\newcommand{\een}{\end{enumerate}}

\title{Amalgamated Products of Groups II: Measures of Random Normal Forms}

\author{\textsf{Elizaveta Frenkel}\\ \textsf{Alexei G. Myasnikov}\\ \textsf{Vladimir N. Remeslennikov\thanks{Supported by
RFFI grant 08-01-00067-a.}}}

\begin{document}
\maketitle

\begin{abstract}
Let $G=\mathop{A\ast B}\limits_C$ be an amalgamated product of
finite rank free groups $A$, $B$ and $C$. We introduce atomic
measures and corresponding asymptotic densities on a set of normal
forms of elements in $G$.  We also define two strata of normal
forms: the first one consists of regular (or stable) normal forms,
and second stratum is formed by singular (or unstable) normal
forms. In a series of previous work about classical
algorithmic problems, it was shown that standard algorithms work
fast on elements of the first stratum and nothing is known about their
work on the second stratum. In theorems \ref{nf:basic} and
\ref{nf_lambdaLmeasurable} of this paper we give probabilistic and
asymptotic estimates of these strata.

\end{abstract}

\begin{minipage}{5in}
\tableofcontents
\end{minipage}


{\bf {Introduction.}} Let $F$ be a free group with basis $X$, $|X|
< \infty$. In \cite{multiplicative,fmrI} the authors introduced a
technique which help to analyze complexity of algorithmic problems
for finitely generated groups of type $G=F/N$. In practical
computations elements of $G$ are usually written in a form of
freely-reduced words in $X$ (normal forms in $G$), and therefore
all computations in $G$ take place in $F$. To analyze a given
algorithmic problem in $G$ one should have:

\begin{itemize}
\item[(i)] satisfactory normal forms for elements of $G$;

\item[(ii)] convenient generators of random elements of $G$ in
normal forms;

\item[(iii)] atomic and probability measures on $F$ for measuring
elements and subsets of $F$;

\item[(iv)] results on stratification of inputs of algorithms
(i.e. normal forms) at least on two strata: stratum of regular (or
stable) elements on which algorithms work fast, for example, in
polynomial time, and another one of singular and unstable elements
on which the result of the algorithms work is unknown or it works
slow;

\item[(v)] asymptotic and probabilistic tools of estimation of
these strata.

\end{itemize}

Here we lay out briefly, what has been done in previous papers of the
authors (and their coauthors) and main results of this paper.

First of all, we work with groups representable in a form of some
free construction, mostly as a free product with amalgamation,
i.e. $G = \mathop{A\ast B}\limits_{C}$. It guarantees the existence of
convenient normal forms of elements in $G$ if these forms exist in
$A$ and $B$. If $A, B, C$ are free groups of finite ranks, we
construct (see section \ref{Section:generators} for details) four
generators of random elements in different normal forms and
specify probabilities to obtain such elements.

In \cite{multiplicative} there was constructed a series of atomic
measures $\{\mu_s| 0 < s < 1, s \in \mathbb{R}^+\}$ on $F$ with
the help of a no-return random walk on the Cayley graph of
$F=F(X)$. It allows us, firstly, to make an asymptotic
classification of subsets of $F$, and, secondly, to prove an
important result about asymptotic properties of regular subsets of
$F$ (i.e. sets accepted by finite automaton).

{\bf {Theorem 1.}}\cite[Theorem
3.2]{multiplicative}.\label{th:thick} { \it Let\/ $R$ be a regular subset of\/
$F$. Then $R$ is thick if and only if its prefix closure
$\overline{R}$ contains a cone{\footnote{All necessary definitions
we give below in Section \ref{section:estimates}.}}.}

In \cite{fmrI} this result was generalized to a stronger form:

{\bf {Theorem 2.}} \cite[Theorem 5.4]{fmrI}. {\it Let $R$ be a
regular subset of a prefix-closed regular set $L$ in a finite rank
free group $F$. Then either the  prefix closure $\overline{R}$ of
$R$ in $L$ contains a non-small  $L-$cone or $\overline{R}$ is
exponentially $\lambda_L$-measurable.}

This theorem plays a significant role in the proof of the main results
(Theorems \ref{nf:basic} and \ref{nf_lambdaLmeasurable}) of this
paper.

Stratification of inputs was described in the papers \cite{CPAI,fmrI}. In particular, we estimate the sizes of strata in
Schreier systems of representatives (transversals) in a free
group. The main result here is a following

{\bf{Theorem 3.}} \cite[Theorem 5.4]{fmrI}. {\it Let $C$ be a
finitely generated subgroup of infinite index in $F(X)$ and $S$ be
a Schreier transversal for $C$. Then sets of all singular
representatives $S_{\rm sin}$ and unstable representatives $S_{\rm
uns}$ are exponentially negligible relative to $S$.}

An example of a stratification of inputs is an algorithm deciding the
Conjugacy Search problem for elements of a group given in
\cite{CPAI}:

{\bf{Theorem 4.}} \cite[Corollary
4.19]{CPAI}.\label{cpaI_cor4.19} {\it Let $G = \mathop{A\ast
B}\limits_{C}$ be a free product of finitely generated free groups
$A$ and $B$ with amalgamated finitely generated subgroup $C.$ Then
the Conjugacy Search Problem in $G$ is decidable for all
cyclically reduced canonical forms $G.$}

Here the cyclically reduced regular elements form the first stratum
and complementary set form the second one.

Main results of this paper is related to item (v) of the research
program described above and is contained in the two theorems:

{\bf{Theorem $A$. }} {\it Let $G=\mathop{A\ast B}\limits_C$ be an
amalgamated product, where $A,B,C$ are free groups of finite rank.
Then for every set of normal forms $\mathcal{NF} = \{
\mathcal{EF}, \mathcal{RF}, \mathcal{CNF}, \mathcal{CRF} \}$

\bi

\item [(i)] If $C$ has a finite index in $A$ and in $B,$ then
every normal form is singular and unstable, i.e.
$\mathcal{NF}_{\rm sin} = \mathcal{NF}_{\rm uns} = \mathcal{NF};$

\item [(ii)] If $C$ of infinite index either in $A$ or in $B,$
then $\mathcal{NF}_r$ and $\mathcal{NF}_s$ are exponentially
$\mu-$generic relative to $\mathcal{NF},$ and $\mathcal{NF}_{\rm
sin}$ and $\mathcal{NF}_{\rm uns}$ are exponentially
$\mu-$negligible relative to $\mathcal{NF}$ in the following
cases:

\bi \item [(ii.1)] $\mu$ is defined by pseudo-measures $\mu_A$ and
$\mu_B$, which are cardinality functions on $A$ and  $B$
correspondingly; in this case $\rho_{\mu}$ is a bidimensional
asymptotic density;

\item [(ii.2)] $\mu$ is defined by atomic probability measures
$\mu_{A,l}$ and $\mu_{B,l}$ on $A$ and $B$ correspondingly; in
this case $\rho^C$ is a bidimensional Cesaro asymptotic density.
\ei \ei }

{\bf {Theorem $B$. }}{\it Let $G=\mathop{A\ast B}\limits_C$ be an
amalgamated product, where $A,B,C$ are free groups of finite rank.
If $C$ of infinite index either in $A$ or in $B,$ then sets of all
unstable $\mathcal{NF}_{\rm uns}$ and all singular
$\mathcal{NF}_{\rm sin}$ normal forms are exponentially
$\lambda_{\mathcal{NF}}-$measurable, where $\mathcal{NF} = \{
\mathcal{EF}, \mathcal{RF}, \mathcal{CNF}, \mathcal{CRF} \}$. }

Specifically, the current papers relation to the other work of this
series is the following. Here, we work with a group $G=\mathop{A\ast
B}\limits_C$ where $A=F(X)$, $B=F(Y)$, and $C$ are free groups of
finite ranks and given atomic measures $\mu_A$ and $\mu_B$ on $A$
and $B$ correspondingly, as well as asymptotic densities induced
by these measures. It was necessary to define correctly
bidimensional measure $\mu$ on $F=A \ast B$ and asymptotic
densities of subsets of $A \ast B$. We do it in Section
\ref{section:asdensities} of this paper.

We follow \cite{LS,MKS} on the subject of group theory;
\cite{eps} for formal languages; \cite{KS,woess} for random
walks. We essentially use the terminology of the papers
\cite{multiplicative,CPAI,fmrI}.

\section{Preliminaries}\label{Section:preliminaries}

In this section we recap some of the definitions and facts about
free products with amalgamation. We refer to \cite{MKS} for more details.
Let $A, B, C$ be groups and $\varphi: C \rightarrow A$ and $\psi: C
\rightarrow B$ be monomorphisms. Then one can define a group $G =
\mathop{A\ast B}\limits_{C}$, called the amalgamated product of $A$
and $B$ over $C$ (the monomorphisms $\varphi, \psi$ are usually
suppressed from notation). If $A$ and $B$ are given by presentations
$A = \langle X \, | \,R_A = 1\rangle,$ $B = \langle Y \,| R_B \, = 1\rangle,$ and a generating
set $Z$ is given for the group $C,$ then the group $G$ has a
presentation

\begin{equation}\label{presentation_G}
G = \langle X  \cup Y | R_A = 1,R_B = 1, \varphi(z) = \psi(z), z \in Z\rangle.
\end{equation}

If we denote $\varphi(z) = u_z(x),\,\, \psi(z)  = v_z(y)$ then $G$
has a presentation
$$G = \langle X \cup Y | R_A = 1,R_B = 1, u_z(x) = v_z(y), (z \in Z )\rangle.$$
Groups $A$ and $B$ are called factors of the amalgamated product $G
= \mathop{A\ast B}\limits_{C};$ they are isomorphic to subgroups in
$G$ generated respectively by $X$ and $Y$. We will identify $A$
and $B$ with these subgroups via evident maps.

Denote by $S$ and $T$ the fixed systems of right coset representatives
of $C$ in $A$ and $B$ respectively. Throughout this paper we assume that the
representative of $C$ is the identity element $1$. For an element $g
\in (A \cup B) \setminus C$ we define $F(g) = A$ if $g \in A$ and
$F(g) = B$ if $g \in B.$

Here are the four main methods to represent an element of $G:$

\begin{enumerate}

\item by a word in the alphabet $X \cup X^{-1}\cup Y
\cup Y^{-1}:$

According to the presentation (\ref{presentation_G}) of $G$, every
nontrivial element $g \in G$ can be written in the form

\begin{equation}\label{any_g}
g = g_1 g_2 \ldots g_n,
\end{equation}
where $g_1, \ldots , g_n$ are reduced words in $X \cup X^{-1}$ or
in $Y \cup Y^{-1}$, and if $F(g_i) = A$ then $F(g_{i+1})=B$, and
vice versa. Obviously, the form (\ref{any_g}) of an element $g$ is not
unique; moreover, the number $n$ of multipliers corresponding to
different representations of $g$ in the form (\ref{any_g}) can
vary ad libitum.

\item in the unique canonical normal form (see \cite{MKS} for
details):
\begin{equation}\label{with}
g = cp_1p_2 \ldots p_l,
\end{equation}
where $c \in C, \,\, p_i \in (S \cup T)\ \setminus \{1\}$, and
$F(p_i) \neq F(p_{i+1}), \,\,i = 1, \ldots, l, \,\, l \geq 0$.

If the {\bf Coset Representative Search Problem} (see, for example,
\cite{CPAI} for details about the algorithmic problems in groups) is
decidable for $C$ in $A$ and $B$ then (\ref{with}) can be computed
from (\ref{any_g}) effectively.

\item in the reduced form:

\begin{equation}\label{red_nf}
g = c g_1 g_2 \ldots g_k
\end{equation}

where $c \in C,$ $g_i \in (A \cup B) \setminus C$ and $F(g_i) \neq
F(g_{i+1}), \,\, i = 1, \ldots, k,$ if $k \geq 0.$ This form may
not be unique, but the number $k$ is uniquely determined by $g.$
For technical reasons, we will use a slightly different definition of
a reduced form as well. Namely, the element $g \in G$ is written in a
reduced form, if
\begin{equation}\label{red_nf_without}
g = g_1 g_2 \ldots g_k
\end{equation}

where $g_i \in (A \cup B) \setminus C$ and $F(g_i) \neq
F(g_{i+1}), \,\, i = 1, \ldots, k,$ if $k \geq 1$ and $g = c$, if
$k = 0.$

Obviously, both definitions (\ref{red_nf}) and
(\ref{red_nf_without}) are equivalent. Moreover, if the {\bf
Membership Problem} for $C$ in $A$ and $B$ is decidable, then
(\ref{red_nf_without}) can be computed from (\ref{any_g})
effectively. Every element in the reduced form (\ref{red_nf}) is a
conjugate of an element

\item in the cyclically reduced form:

\begin{equation}\label{crf} g=c g_1 \ldots g_k
\end{equation}
The form (\ref{crf}) is called \emph{cyclically reduced form} of
element $g,$ if

\bi

\item [(i)] $k = 0,$ i.e. $g=c \in C;$

\item [(ii)] $k = 1,$ then every $g \in A \cup B$ that
 is not a conjugate of an element in $C;$

\item [(iii)] $k > 1,$ then every $g,$ such as $k$ is even.

\ei


\end{enumerate}
We will refer to an element $g \in G$ as an element in normal
form throughout the paper if it has one of the forms (\ref{any_g}),
(\ref{with}), (\ref{red_nf_without}) or (\ref{crf}) and it doesn't
matter which one is chosen. Though, we mention that the numbers $n, k,
l$ of the factors in each representation of an element in different normal forms
are different in general (and some time are not even unique), we will
refer further to such a number as to length of the representation of $g$
in normal form and denote it by $s(g)$.

\subsection{Measuring and comparing subsets of free
group}\label{subsection:measuringF}

In this section we recap some crucial facts about measures in
free group of finite rank $F(X);$ in more details you can find
this information in \cite{multiplicative,fmrI}.

Let $\mathcal{P}(F)$ be the set of all subsets of $F = F(X)$ and
$\mathcal{A} \subset \mathcal{P}(F).$ A real-valued non-negative
additive function $\mu: \cal{A} \rightarrow \mathbb{R}^{+}$ is
called a \emph{pseudo-measure on $F$}. If $\mathcal{A}$ is a
subalgebra of $\mathcal{P}(F)$, then $\mu$ is a \emph{measure}.

Let $F=F(X)$ be a free group. Denote by $S_n$ and $B_n$
correspondingly the sphere and the ball of radius $n$ in $F.$ Let
$\mu$ be an atomic pseudo-measure on $F.$ Recall, that a measure
$\mu$ on the countable set $P$ called is \emph{atomic} if every
subset $Q \subseteq P$ is measurable; it also holds when
$\mu(Q)=\sum\limits_{q\in Q}\mu(q).$

For a set $R \subseteq F$ we define its \emph{spherical asymptotic
density} relative to $\mu$ as the following limit:

$$s\rho_{\mu}(R) = \overline{\lim\limits_{n \rightarrow \infty}} s\rho_n(R),$$
where $s\rho_n(R) =  \fracd{\mu (R \bigcap S_n)}{\mu(S_n)}.$

Similarly, one can define the \emph{ball asymptotic density} of
$R$ relative to $\mu$:
$$b\rho_{\mu}(R) = \overline{\lim\limits_{n \rightarrow \infty}} b\rho_n(R),$$
where $b \rho_n(R) = \fracd{\mu(R \bigcap B_n)}{\mu(B_n)}.$ We
formulate below one very useful fact about the connection between
spherical and ball asymptotic densities, a proof can be found, for
example, in \cite[Lemma 3.2]{GMMU}.
\begin{lemma}\label{lemma1}
Let $\mu$ be a pseudo-measure on $F.$ Suppose that
$\mathop{lim}\limits_{n \rightarrow \infty} \mu(B_n) = \infty.$

Then for any subset $R\subseteq F$ if the spherical asymptotic
density $s\rho_{\mu}(R)$ exists, then the ball asymptotic density
$b\rho_{\mu}(R)$ also exists and
$$b\rho_{\mu}(R)=s\rho_{\mu}(R).$$
\end{lemma}

Further, let $\mu$ be a pseudo-measure on $F$ and $\rho_{\mu}(R)$ be
a spherical or ball asymptotic density.

We say that a subset $R \subseteq F$ is \emph{generic } relative
to $\mu$, if the limit $\lim\limits_{n \rightarrow \infty}
s\rho_n(R)$ exists and $\rho_{\mu}(R)=1$, and \emph {negligible }
relative to $\mu$, if $\rho_{\mu}(R)=0$.

Further, we say $R$ is \emph{exponentially generic} relative to
$\mu$ if there exists a positive constant $\delta <1$ such that
$1-\delta^n < s\rho_{n}(R) < 1$ for sufficiently large $n$.
Meanwhile, if $s\rho_{n}(R) < \delta^n$ for all large enough $n$,
then $R$ is {\em exponentially negligible} relative to $\mu$.

For example, if $\mu$ is the cardinality function, i.e.
$\mu(A)=\abs{A},$ then we obtain standard asymptotic density
functions on $F$. We will use the notation $\rho(R)$ throughout
the paper to denote the standard spherical asymptotic density of
$R$ in $F$ relative to the cardinality function; we will also omit the
''cardinality function'' whenever possible. It is not also hard
to extend such a definition for an asymptotic density of set $R$
relative to set $R_1$ (see \cite{fmrI} and Section
\ref{subsection:lamda_l} for details). We will use the notation
$\rho_{\mu}(R,R_1)$ for this asymptotic density.

Further, we will be interested on a special kind of measure in $F,$
studied in details in a lot of papers (see, for example,
 \cite{af_semr,multiplicative,fmrI}). Namely, consider a so-called {\em
frequency measure} on $R \subseteq F = F(X):$
$$\lambda(R) = \sum_{n=0}^{\infty}f_n(R), \textrm{ where } f_n(R)= \frac{\vert R \cap S_n\vert}{\vert S_n\vert},$$
 and $f_n(R)$ are called {\em  frequencies} of elements from $R$ among the words of
(freely-reduced) length $n$ in $F.$ This measure  is not
probabilistic, since, for instance, $\lambda(F) = \infty$,
moreover, $\lambda$ is additive, but not $\sigma$-additive.

 Also, frequencies of $R$ define a well-studied
asymptotic density called {\em Cesaro asymptotic density}. Namely,
it is the {\em Cesaro limit} of frequencies for $R:$
 \begin{equation}
 \label{eq:cesaro-1}
  \rho^c(R)  = \lim_{n\rightarrow\infty} \frac{1}{n}\left(f_1(R)+\cdots+f_n(R)\right).
  \end{equation}
 Sometimes it is more sensitive then the standard asymptotic density $\rho$ (see,
for example, \cite{multiplicative,woess}). However, if
$\lim\limits_{n\rightarrow \infty} f_n(R)$ exists (hence is equal
to $\rho(R)$), then $\rho^c(R)$ also exists and $\rho^c(R) =
\rho(R).$

\subsection{Stratification and measuring of Schreier systems of representatives in free group}\label{subsection:measuringS}
In this section we give some information about Schreier
transversals (see \cite{fmrI} for details) in free groups and also
the definitions of regular and stable normal forms of elements in
free product with amalgamation $G.$

Following  \cite{km}, we associate with $C$ two graphs: the
\emph{subgroup graph} $\Gamma =\Gamma_C$ and the Schreier graph
$\Gamma^{\ast} = \Gamma^{\ast}_C$. Recall that $\Gamma$  is a
finite connected digraph with edges labeled by elements from $X$
and a distinguished vertex (based-point) $1_C$,  satisfying  the
following two conditions. Firstly,  $\Gamma$ is folded, i.e.,
there are no two edges in $\Gamma$ with the same label and having
the same initial or terminal vertices. Secondly, $\Gamma$ accepts
precisely the reduced words in $X \cup X^{-1}$  that belong to
$C$.

The {\em  Schreier graph} $\Gamma^{\ast} = \Gamma_C^{\ast}$ of $C$
is a connected labeled digraph with the set $\{Cu \mid u \in F\}$
of right cosets of $C$ in $F$ as the vertex set, and such that
there is an edge from
 $Cu$ to $Cv$ with a label $x \in X$ if and only if $Cux = Cv$.  One can describe the Schreier graph $\Gamma^{\ast}$ as  obtained from $\Gamma$
by the following procedure. Let $v \in \Gamma$ and $x \in X$ such
that there is no outgoing or incoming
 edge at $v$ labeled by $x.$ For every such vertex $v$ and $x \in X$ we attach
to $v$ a new edge $e$ (correspondingly, either outgoing or
incoming) labeled $x$ with a new terminal vertex $u$ (not in
$\Gamma$). Then we attach to $u$ the Cayley graph $C(F,X)$ of $F$
relative to $X$ (identifying $u$ with the root vertex of
$C(F,X)$), and
 then we fold the edge $e$ with the corresponding edge in $C(F,X)$ (that is labeled $x$ and is incoming to $u$).
Observe, that for every vertex $v \in \Gamma^{\ast}$ and every
reduced word $w$ in $X \cup X^{-1}$ there is a unique path
$\Gamma^{\ast}$ that starts at $v$ and has the label $w$. By $p_w$
we denote such a path that starts at $1_C$, and by $v_w$ the end
vertex of $p_w$.

Consider a set of right representatives of $C$ in $F= F(X)$; we
will call it the {\em transversal} of $C.$ Recall, that a
transversal $S$ of $C$ is termed {\em Schreier} if  every initial
segment of a representative from $S$ belongs to $S$. In
\cite{fmrI} was shown, that there is one-to-one correspondence
between the set of every Schreier transversal $S$ of $C$ and the
set of all spanning subtrees $\Gamma^{\ast}.$ In particular, it
means that we can treat with every representative $s \in S$ as
with label of a path in some (fixed) spanning subtree of
$\Gamma^{\ast}.$

Also, we have a classification of representatives of $C$ in $F=F(X)$
from \cite{fmrI}:

\begin{definition} Let $S$ be a transversal of $C$.
\begin{itemize}

\item
 A representative $s\in S$ is
called \emph{internal} if the path $p_s$ ends in $\Gamma$, i.e.,
$v_s \in V(\Gamma)$.
 By $S_{\rm int}$ we denote the
set of all internal representatives in $S.$  Elements from $S_{\rm
ext}=S \smallsetminus S_{\rm int}$ are called \emph{external}
representatives in $S.$

\item A representative $s \in S$ is called {\em singular} if it belongs to the generalized normalizer of $C$:
$$N^*_F(C) = \{ f \in F | f^{-1} C f \cap C \neq 1 \}.$$
 All other representatives
from $S$ are called {\em regular}. By $S_{\rm sin}$ and,
respectively, $S_{\rm reg}$ we denote the sets of singular and
regular representatives from $S$.

\item A representative $s \in S$ is called {\em stable} if $sc \in
S$ for any $c \in C$.  By $S_{\rm st}$ we denote the set of all
stable representatives in $S,$ and $S_{\rm uns} = S \smallsetminus
S_{\rm st}$ is the set of all {\em unstable } representatives from
$S$.
\end{itemize}
\end{definition}

\emph{Frontier} vertex $v_u$ of $V(\Gamma)$ is a vertex $v_u \in
V(\Gamma^{\ast}) \setminus V(\Gamma)$ such that $v_u$ incident to
an edge $e$ of $\Gamma^{\ast}$, which initial or terminal vertex
already in $V(\Gamma)$. A \emph{cone} $C(u)$ is a subset of $F$ of
type $\{ w \in F: w = u f {\textrm{ and }} uv {\textrm{ is a
reduced word } }\}$.

The following proposition about the structure of all singular and
unstable representatives was shown in \cite{fmrI}:
\begin{proposition}
\label{pr:basic-properties}\cite[Proposition 3.5]{fmrI}. Let $S$
be a Schreier transversal for $C$, $C$ has an infinite index in $F$
and $S = S_{T^{\ast}}$ for some spanning subtree ${T^{\ast}}$ of
$\Gamma^{\ast}$. Then the following hold:
  \bi
   \item [1)] $|S_{\rm int}| = |V(\Gamma)|.$
    \item [2)] $S_{\rm ext}$ is  the union of finitely many coni $C(u),$ where $v_u$ are frontier vertices of $\Gamma$.

    \item  [3)]$S_{\rm sin}$ is contained in a finite union of double
    cosets $Cs_1s_2^{-1}C$ of $C$, where $s_1, s_2 \in S_{\rm int}$.
    \item [4)] $S_{\rm uns}$ is a finite union of  left
    cosets of $C$ of the type $s_1s_2^{-1}C$, where $s_1, s_2 \in S_{\rm
    int}$.
    \item [5)] (see [Proposition 3.9, \cite{fmrI}]) $S_{\rm sin} \subseteq S_{\rm uns}$.

 \ei
\end{proposition}

 In \cite{fmrI} it was also shown that the sets of all
singular and unstable representatives forms an exponentially
negligible part of the Schreier transversals:
\begin{corollary}\label{s_nst-ins}\cite[Corollary 5.12]{fmrI}.
Let $C$ be a finitely generated subgroup of infinite index in
$F(X)$ and $S$ a Schreier transversal for $C$. Then sets of
singular representatives $S_{\rm sin}$ and unstable
representatives $S_{\rm uns}$ are exponentially negligible in $S$.
\end{corollary}

Following the idea to split the set of all normal forms of
elements in $G = \mathop{A\ast B}\limits_C$ into "bad" and "good"
components, we
introduce the following definitions. 

We say that an element $g \in F(X)$ is regular (stable and so on),
if it can be decomposed into a form $g = c s, \,\, c \in C, s \in S$
and $s$ is regular, stable etc.

\begin{definition} An element $g \in
G$ in normal form (\ref{any_g}), (\ref{with}) ,
(\ref{red_nf_without}) or (\ref{crf}) is called \emph{regular} if at
least one of elements $g_i \textrm{ or }\,\, p_i; \,\, i = 1,
\ldots, s(g)$ is regular. Otherwise $g$ is called \emph{singular}.
\end{definition}

\begin{definition}
An element $g \in G$ in normal form (\ref{any_g}), (\ref{with}),
(\ref{red_nf_without}) or (\ref{crf}) is called \emph{stable} if at
least one of elements $g_i \textrm{ or }\,\, p_i; \,\, i = 1,
\ldots, s(g)$ is stable. Otherwise $g$ is called \emph{unstable}.
\end{definition}

In main the Theorems \ref{nf:basic} and \ref{nf_lambdaLmeasurable} of
this paper we estimate sizes of stable, unstable, regular and
singular components in the set of all normal forms, and these
notions will be very important for the rest of the paper.

\section{Asymptotic densities on free products of
subsets}\label{section:asdensities}

In this Section our goal is a definition of asymptotic densities
on subsets of $G= A \mathop{\ast}\limits_C B$, induced by
different types of measures on factors $A$ and $B$, introduced in
Section \ref{subsection:measuringF}. In turn, it constrain us to
define measures on free product of subsets of $F = A \ast B.$
\subsection{Free product of subsets}
Let $F= A \ast B$ be a free product of finitely generated groups
$A$ and $B$, and $A_0$ is a nonempty subset of $A,$ $B_0$ is a
nonempty subset of $B$. Then a \emph{free product} $A_0 \ast B_0$
is a set of all elements in $F$ having a form $f= f_1 f_2 \ldots
f_k,$ where $n \geq 1; \,\,f_i \in A_0 \cup B_0, \,\, i = 1,
\ldots, k,\,\, f_i \neq 1$ if $i \geq 2$ and for all $i = 1,
\ldots, k-1$ elements $F(f_i) \neq F(f_{i+1})$. Every nontrivial
element $f \in F$ can be written in a freely-reduced form
\begin{equation}\label{any_f}
f = f_1 f_2 \ldots f_k,
\end{equation}
where $f_1, \ldots , f_n$ are reduced words in $X \cup X^{-1}$ or
in $Y \cup Y^{-1}$, and if $F(f_i) = A$, then $F(f_{i+1})=B$ and
vice versa. For every such nontrivial $f \in F$ set $s(f)=k$; let
$s(1) =0$. We will denote by $\abs{f},\,\abs{f_i}$ the number of
letters in alphabet $X \cup X^{-1}$ or $X \cup X^{-1} \cup Y \cup
Y^{-1}$ in a freely-reduced form of words $f, f_i$.

Let $\mu_A$ and $\mu_B$ be atomic pseudo-measures on $A$ and $B$
correspondingly and $\mu_A(1) = \mu_B(1).$ We also fix a
probability distribution $\theta: \mathbb{N} \rightarrow
\mathbb{R}^+,$ in particular, $\mathop{\sum}\limits_{k
=1}^{\infty}\theta(k) = 1.$ We define an atomic measure $\mu$ on
$F$ in the following manner:
\begin{equation}\label{mu_f}
\mu(f)= \frac{1}{2} \theta (k) \mu_{F_1}(f_1) \ldots
\mu_{F_k}(f_k),
\end{equation}
where $f$ is written in a freely-reduced form (\ref{any_f}).

For a subset $R \subseteq F$ set
$$\mu(R) = \mathop{\sum}\limits_{f \in R} \mu(f).$$
We will say, that $R$ is a \emph{$\mu-$measurable} set, if $\mu(R) <
\infty.$ Denote by $M_{\mu}$ the set of all $\mu-$measurable subsets
of $F:$
$$M_{\varphi} = \{ R
\subseteq F | \mu(R) < \infty\}.$$

\begin{lemma}\label{lem:free_product}Let $F = A \ast B$ be a free product of finitely generated groups $A$
and $B$.
 \bi \item [1)] If $\mu_A, \mu_B$ are atomic
pseudo-measures on $A$ and $B$ correspondingly, then measure $\mu$
on $F$ defined above is an atomic pseudo-measure on $F$ and
$M_{\mu_A} \subset M_{\mu}, M_{\mu_B} \subset M_{\mu}.$

\item [2)] If $\mu_A, \mu_B$ are atomic probability measures on
$A$ and $B$ correspondingly, then $\mu$ is an atomic probability
measure on $F.$ \ei
\end{lemma}
\proof Claim 1) is straightforward. Let us prove 2).
Split $F$ into layers: $F = F_1 \bigsqcup F_2 \bigsqcup F_3 \ldots
,$ where
$$F_i = \{ f \in F | f \textrm{ in a freely-reduced form (\ref{any_f}) and }
s(f) = i \}.$$
$$\begin{array}{l}
\textrm{Then }\,\,\, \mu(F) = \mathop{\sum}\limits_{i = 1}^{\infty}\mu(F_i) =\\
=\mathop{\sum}\limits_{i = 1}^{\infty} \fracd{1}{2}
\theta(i)\left(\mu_A^{[\frac{i+1}{2}]}(A)\mu_B^{i-[\frac{i+1}{2}]}(B)+
\mu_A^{i-[\frac{i+1}{2}]}(A)\mu_B^{[\frac{i+1}{2}]}(B)\right) =\\
=\mathop{\sum}\limits_{i = 1}^{\infty} \theta(i) = 1.
\end{array}$$ \eproof

\textbf{Example.} Suppose $\mu_A, \mu_B$ are pseudo-measures on
$A$ and $B$, defined by cardinality functions on $A$ and $B$, and
$\theta(k)= \fracd{6}{\pi^2 k^2}$ is a probability distribution on
$\mathbb{N.}$ Then $M_{\mu_A} = \mathcal{F}(A)$ and $M_{\mu_B} =
\mathcal{F}(B),$ where $\mathcal{F}(A)$ and $\mathcal{F}(B)$ are
sets of all finite subsets of $A$ and $B.$ However, $M_{\mu}
\supset \mathcal{F}(F)$ is a strict inclusion. Indeed, let $R
\subseteq F$ and $R_k = R \cap F_k.$ Then $R \in M_{\mu}$ iff row
$\mathop{\sum} \fracd{\abs{R_k}}{k^2}$ converges.

We shall describe below several methods to define asymptotic
density of subsets in $F = A \ast B$.

\subsection{Bidimensional asymptotic density}
Let $T= A_0 \ast B_0 \subseteq F=A \ast B.$ For a pair of natural
numbers $(n,k)$ we define \emph{$(n,k)$-ball}:
$$T_{n, k}=\{ f = f_1 \ldots f_k \in T :  s(f) = k, \abs{f_i} \leq n, i = 1,\ldots, k \}.$$
We call $T= \mathop{\cup}\limits_{k
=0,n=0}^{\infty}T_{n,k}$ the
 \emph{bidimensional decomposition} of $T$.
Bidimensional decompositions help us to analyze asymptotic
behavior of subsets of $T$ and other subsets of $F$ relative to
$T$. For a set $Q=A_1 \ast B_1$ in $F$ a function $(n,k)
\rightarrow \mu(Q \cap T_{n,k})$ is called the growth function of
$Q$ in $T,$ and a function $(n,k) \rightarrow
\rho_{\mu}^{n,k}(Q,T) = \fracd{\mu (Q \cap
T_{n,k})}{\mu(T_{n,k})}$ is called the frequency function of $Q$
relative to $T.$

By \emph{direction function} $d(n,k)$ we mean one-to-one
correspondence between $n$ and $k$ which parametrize a path from
$(1,1)$ to $(\infty, \infty)$ such that arguments $n$ and $k$
tends to $\infty$ while $d(n,k) \rightarrow \infty$.
 Let $d(n,k)$ be some direction function.
Asymptotic behavior of $Q$ relative to $T$ we will characterize by
a \emph{bidimensional asymptotic density}, which determines as
following limit:
$$\rho_{\mu}(Q,T) = \mathop{\overline{\lim}}\limits_{d(n,k) \rightarrow \infty} \rho_{\mu}^{n,k}(Q,T).$$
If this limit exists and does not depend on a choice of a
direction function, we denote it by $\rho^e_{\mu}(Q,T)$. We say
that $Q$ is \emph{$\mu$-generic} relative to $T$ , if
$\rho^e_{\mu}(Q,T)=1$, and \emph{$\mu$-negligible} relative to
$T$, if $\rho^e_{\mu}(Q,T)=0$.

Further, we say that $Q$ is \emph{exponentially $\mu_n-$generic}
relative to $T$ if $Q$ is $\mu-$generic and there exists a
positive constant $\delta <1$ such that $1-\delta^{nk} <
\rho_{\mu}^{n,k}(Q,T) < 1$ for all large enough $n$ and $k$.
Meanwhile, if $Q$ is $\mu$-negligible and $\rho_{\mu}^{n,k}(Q,T) <
\delta^{nk}$ for sufficiently large $n, k$ then $Q$ is {\em
exponentially $\mu_{n}-$negligible} relative to $T$.

We will use below other notions of exponentially generic and
negligible sets relative to measure $\mu$, that will allow us to
obtain in Section \ref{section:estimates} more rough, but at the
same time, more general estimates on subsets of normal forms.
Namely, $Q$ is said to be {\em exponentially
 $\mu-$generic} relative to $T$ if $Q$ is $\mu-$generic and there
exists a positive constant $\delta <1$ such that $1-\delta^{k} <
\rho_{\mu}^{n,k}(Q,T) < 1$ for large enough $k$. If $Q$ is
$\mu$-negligible and $\rho_{\mu}^{n,k}(Q,T) < \delta^{k}$ for
sufficiently large $k$, then $Q$ is {\em exponentially
$\mu-$negligible} relative to $T$.

Denote $(T)_n = \{f = f_1\ldots f_k \in T \,\,|\,\,\, \abs{f_i} =
n,\,\,i=1,\ldots,k\}$,  $(T)_{\leq n} = \{f = f_1\ldots f_k \in T
\,\,|\,\,\, \abs{f_i} \leq n,\,\,i=1,\ldots,k\}$,  and $(T)^k =
\{f \in T \,\,|\,\,\, s(f) = k \}$. We shall use the same notation
for subsets of groups $A$ and $B$ when there is no ambiguity, i.e.
notation $(A_i)_n = (A)_n \cap A_i$ (or $(B_i)_n = (B)_n \cap
B_i$) for the sphere of radius $n$ in a subsets of $A$ (or $B$)
and $(A_i)_{\leq n} = (A)_{\leq n} \cap A_i$ (or $(B_i)_{\leq n} =
(B)_{\leq n} \cap B_i$) for the ball of radius $n$ in
corresponding sets.

The following proposition will be useful in the sequel.
\begin{proposition}\label{munegligibility}
Let $F=A \ast B$ be a free product of free groups $A$ and $B$ of
finite ranks, and $A_1 \subseteq A_0 \subseteq A,$ $B_1 \subseteq
B_0 \subseteq B$, and let $T=A_0 \ast B_0$, and $Q=A_1 \ast B_1$.
\begin{itemize}
\item[1.] Suppose $\mu_A$ and $\mu_B$ are atomic probability
measures on $A$ and $B$ correspondingly and $\mu$ is a
pseudo-measure on $F$, defined above and let
$\rho_{\mu_A}(A_1,A_0)$ and $\rho_{\mu_B}(B_1,B_0)$ exist. If
there is a constant $\delta, 0< \delta <1$ such that for all $n >
n_0$ either $\fracd{\mu_A ((A_1)_{\leq n})}{\mu_A ((A_0)_{\leq
n})} < \delta^n$ or $\fracd{\mu_B ((B_1)_{\leq n})}{\mu_B
((B_0)_{\leq n})} < \delta^n$, then $Q$ is exponentially
$\mu_{n}-$negligible relative to $T$.

\item[2.] Suppose $\mu_A$ and $\mu_B$ are atomic probability
measures on $A$ and $B$ correspondingly and $\mu$ is a
pseudo-measure on $F$, defined above. Let $\rho_{\mu_A}(A_1,A_0)$
and $\rho_{\mu_B}(B_1,B_0)$ exist and at least one of them less
than $1$. Then $Q$ is exponentially $\mu-$negligible relative to
$T$.

\item[3.] Suppose $\mu_A$ and $\mu_B$ are pseudo-measures on $A$
and $B$ defined by cardinality functions and $\mu$ is a
pseudo-measure on $F$, defined above. Let $\rho(A_1,A_0)$ and
$\rho(B_1,B_0)$ exist and at least one of them less than $1$. Then
$Q$ is exponentially $\mu-$negligible relative to $T$.
\end{itemize}

\end{proposition}
\proof To prove first claim, suppose $\fracd{\mu_A ((A_1)_{\leq
n})}{\mu_A ((A_0)_{\leq n})} < \delta^n$; then by definition for
all $n \geq n_0$ we obtain $\fracd{\mu (Q \cap
T_{n,k})}{\mu(T_{n,k})}=$

\begin{equation}\label{eq}
 =\left \{
\begin{array}{l}
\textrm{ если } \,\,k =
2t\,\,\,\,\,\,\,\,\,\fracd{(\mu_A((A_1)_{\leq n}))^t
(\mu_B((B_1)_{\leq n}))^t}{(\mu_A((A_0)_{\leq n}))^t (\mu_B((B_0)_{\leq n}))^t}\\
\\
\textrm{ если } \,\,k = 2t+1
\,\,\,\,\,\,\,\fracd{(\mu_A((A_1)_{\leq n}))^t
(\mu_B^t((B_1)_{\leq n}))^t (\mu_A ((A_1)_{\leq
n})+\mu_B((B_1)_{\leq n}))}{(\mu_A((A_0)_{\leq n}))^t
(\mu_B((B_0)_{\leq n}))^t (\mu_A
((A_0)_{\leq n})+\mu_B((B_0)_{\leq n}))}.\\
\end{array}
\right.\end{equation}

Since $\rho_{\mu_B}(B_1,B_0)$ exists there is a natural number
$n_1$ such that $$\fracd{\mu (Q \cap T_{n,k})}{\mu(T_{n,k})} \leq
\delta^{nt} < \delta^{n(t-1)} < \delta^{n(k/2-2)}$$ for all $n
\geq n_1$.
   Therefore, the limit $\mathop{lim}\limits_{d(n,k) \rightarrow \infty}\fracd{\mu (Q \cap T_{n,k})}{\mu(T_{n,k})}$
   exists and does not depend on a $d(n,k)$. It equals to zero and moreover,
 $Q$ is exponentially $\mu_n-$negligible relative to
$T$ with a $\delta' = \delta^{1/2}$.

To prove 2), suppose that $\rho_{\mu_A}(A_1,A_0) < 1$. By simple
observation since the limit exists
$$\rho_{\mu_A}(A_1,A_0) = \mathop{lim}\limits_{n \rightarrow
\infty} \fracd{\mu_A ((A_1)_n)}{\mu_A ((A_0)_n)}=
\mathop{lim}\limits_{{n \rightarrow \infty}}\fracd{\mu_A
((A_1)_{\leq n})}{\mu_A ((A_0)_{\leq n})}$$ and therefore
$$\fracd{\mu_A ((A_1)_{\leq n})}{\mu_A ((A_0)_{\leq
n})} < 1 - \varepsilon$$  for
 relevant $0 < \varepsilon < 1$.

Again, using (\ref{eq}), obtain  $\fracd{\mu (Q \cap
T_{n,k})}{\mu(T_{n,k})}\leq (1-\varepsilon)^t <
(1-\varepsilon)^{t-1} < (1-\varepsilon)^{k/2-2}$ and therefore $Q$
is $\mu-$negligible relative to $T$ for arbitrary choice of a
direction function $d$. It is clear also that $Q$ is exponentially
$\mu-$negligible relative to $T$ with a $\delta =
(1-\varepsilon)^{1/2}$.

The proof of the last claim is analogous to the former one. \eproof

\subsection{Bidimensional Cesaro asymptotic density}

Let $\{\mu_s,\, 0 < s <1 \}$ be a family of probabilistic
distributions introduced in \cite{multiplicative} for a free group
$F(X)$ of finite rank. In terms of relative frequencies of $R$
relative to $F$ it can be written as follows $$ \mu_s(R) =
s\sum_{k=0}^\infty f_k(1-s)^k.$$ In \cite{multiplicative} it was also
shown that the average (freely-reduced) length of words in $F(X),$
distributed according to $\mu_s$ is equal to $l = \fracd{1}{s}-1$;
evidently, $l \rightarrow \infty$ while $s \rightarrow 0^+$.
Therefore, the family $\{\mu_s \}$ can be parametrized by $l$ :
$\{\mu_l| \,\,1 < l < \infty\}$.

Suppose $\mu_A=\{\mu_{A,l}\}$ and $\mu_B=\{\mu_{B,l}\},$
$1<l<\infty$ are atomic probability measures on free groups $A$
and $B$ correspondingly. For asymptotic estimates of sets it is
sufficient to assume that $l$ runs over natural numbers. Let
$\mu_l$ be the atomic probability measure on $F$ induced by
$\mu_{A,l}$ and $\mu_{B,l}$ and let $Q = A_1 \ast B_1 \subseteq T
= A_0 \ast B_0 \subseteq F = A \ast B.$ For some choice of
direction function $d(l,k)$ consider a function $(l,k) \rightarrow
\rho_{\mu_l}^{l,k}(Q,T)$ of relative frequencies.

\begin{definition}
Let $Q \subseteq T \subseteq F$ as above. The function
$$(l,k) \rightarrow \fracd{\mu_l((F)^k \cap Q)}{\mu_l((F)^k \cap T)} = \rho_{\mu_l}^{l,k}(Q,T)$$ is
called frequency function of $Q$ relative to $T$. The limit (if it
exists and does not depend on a choice of $d(l,k)$)
$$\rho^C(Q,T) = \mathop{lim}\limits_{d(l,k) \rightarrow \infty}
\rho_{\mu_l}^{l,k}(Q,T),$$ is called the \emph{bidimensional
Cesaro asymptotic density} of $Q$ relative to $T.$\end{definition}

If $s$ is not small (for example, $s \geq \fracd{1}{2}$), then
every set containing $1$ or short elements is not small since
$\mu_s(1) = s$. We avoid this since $l$ runs over natural $l$'s
and $l>1$.

We say that $Q$ is \emph{$C-$negligible} relative to $T$ if
$\rho^C(Q,T) = 0$; further, $Q$ is \emph{exponentially
$C-$negligible} relative to $T$ if $\rho_{\mu_l}^{l,k}(Q,T) <
\delta^k$ for some constant $0<\delta<1$ and all sufficiently
large $k$. Supplements to these sets called \emph{$C-$generic} and
\emph{exponentially $C-$generic} correspondingly.

Here we describe some sufficient conditions on free product of
sets to be exponentially negligible with respect to Cesaro
asymptotic density, which we will use in Section
\ref{section:estimates} to evaluate sizes of subsets of normal
forms.

\begin{proposition}\label{Cnegligibility}
Let $F= A\ast B$ be a free product of two free groups $A$ and $B$ of
 finite ranks. Suppose $\mu_A=
\{\mu_{A,l}\}$ and $\mu_B= \{\mu_{B,l}\}$, $1<l<\infty$ are two
families of atomic probability measures on $A$ and $B$
correspondingly and  $\{\mu_l\}$ is induced family of measures on
$F$. Let $A_1 \subseteq A_0 \subseteq A$, $B_1 \subseteq B_0
\subseteq B$ and densities $\rho_{\mu_{A,l}}(A_1,A_0)$,
$\rho_{\mu_{B,l}}(B_1,B_0)$ exist. If there is a number $0 < q <
1$ such that $\fracd{\mu_{A,l}(A_1)}{\mu_{A,l}(A_0)} < q$ or
$\fracd{\mu_{B,l}(B_1)}{\mu_{B,l}(B_0)} < q$ for all $l >l_0$ for
some $l_0$, then the set $Q = A_1 \ast B_1$ is exponentially
$C-$negligible relative to $T = A_0 \ast B_0.$
\end{proposition}
\proof Let us fix a pair of natural numbers $(l,k)$ such that $l
> l_0$. Suppose that $\fracd{\mu_{A,l}(A_1)}{\mu_{A,l}(A_0)} < q$.
Splitting $T$ into layers as in Lemma \ref{lem:free_product},
obtain
$$\rho^{l,k}_{\mu_l}(Q,T) =$$

$$
\left \{
\begin{array}{l}
\textrm{ if } \,\,k =
2t\,\,\,\,\,\,\,\,\,\fracd{(\mu_{A,l}(A_1))^t
(\mu_{B,l}(B_1))^t }{(\mu_{A,l}(A_0))^t (\mu_{B,l}(B_0))^t}\\
\\
\textrm{ if } \,\,k = 2t+1 \,\,\,\,\,\,\,\fracd{(\mu_{A,l}(A_1))^t
(\mu_{B,l}(B_1))^t
(\mu_{A,l}(A_1)+\mu_{B,l}(B_1))}{(\mu_{A,l}(A_0))^{t}
(\mu_{B,l}(B_0))^t(\mu_{A,l}(A_0)+\mu_{B,l}(B_0))}.\\
\end{array}
\right. $$

Therefore,
$$\rho^{k,l}_{\mu_l}(Q,T) \leq \left( \fracd{\mu_{A,l}(A_1)}{\mu_{A,l}(A_0)}\right)^t \leq q^t$$
and since the limit does not depend on a choice of a direction,
$A_1 \ast B_1$ is exponentially (with $\delta =q^{\frac{1}{2}}$)
$C-$negligible relative to $A_0 \ast B_0.$\eproof

The following lemma connects two types of measuring in free groups
and will be very useful in the sequel:

\begin{lemma}\label{snegl_is_Cnegl}
Let $A_1$ be an exponentially negligible subset relative to a
subset $A_0$ of a finitely generated free group $A = F(X)$. Then
for arbitrary $l_0 > 1$ there is a real number $0< q < 1$ such
that for all $l > l_0$ holds
$$\fracd{\mu_{A,l}(A_1)}{\mu_{A,l}(A_0)}<q.$$
\end{lemma}
\proof Since $A_1$ is exponentially negligible relative to $A_0$
there is a natural number $n_0$ and real number $0 < p < 1$ such
that $f_n(A_1,A_0) = \fracd{|(A_1)_n|}{|(A_0)_n|} < p^n$ for all
$n \geq n_0$. In particular, $|(A_1)_n|  < p^{n_0} |(A_0)_n|$ for
all $n \geq n_0$. It is sufficient to show that for every fixed $s
< s_0$ a number $\fracd{\mu_{A,s}(A_1)}{\mu_{A,s}(A_0)}$ is
bounded above by some positive constant $q < 1$.

By definition we have
\begin{equation}
\fracd{\mu_{A,s}((A_1)_{\leq n})}{\mu_{A,s}((A_0)_{\leq n})} =
\fracd{s \mathop{\sum}\limits_{k=0}^{n}|(A_1)|(1-s)^k}{s
\mathop{\sum}\limits_{k=0}^{n}|(A_0)|(1-s)^k} < \fracd{
\mathop{\sum}\limits_{k=0}^{n} p^{n_0}|(A_0)|(1-s)^k}{
\mathop{\sum}\limits_{k=0}^{n}|(A_0)|(1-s)^k} < p^{n_0}
\end{equation}

 for all $n
\geq n_0, s > s_0$. Without loss of generality one can assume that
$n_0 > 1$, in opposite case add or remove arbitrary element of
length $1$ from $A_1$. Passing to a limit in the inequality
$\fracd{\mu_{A,s}((A_1)_{\leq n})}{\mu_{A,s}((A_0)_{\leq n})} <
p^{n_0}$, obtain $\fracd{\mu_{A,s}((A_1))}{\mu_{A,s}((A_0))} \leq
p^{n_0}$ and therefore $\fracd{\mu_{A,s}(A_1)}{\mu_{A,s}(A_0)}$ is
strictly bounded by positive constant $q=p^{n_0-1}$ for every $s <
s_0$. \eproof

\section{Generation of random normal forms}\label{Section:generators}
Let $G = \mathop{A\ast B}\limits_{C}$ be a free product of a free
group $A$ with a finite base $X$ and a free group $B$ with a finite
base $Y,$ amalgamated over a finitely generated subgroup $C$.

In the series of previous papers \cite{multiplicative, CPAI} we
described some algorithms in free groups of finite rank and
amalgamated products of such groups. The special role in analysis
of computational complexity of algorithmic problem in group $G$
play regular and stable normal forms. So, the goal of this section
is a construction of procedures for generating of random reduced,
random canonical normal and random cyclically reduced forms. The
second goal is asymptotic estimation of sets of all regular,
stable, singular and unstable forms relative to the set of all
forms with the help of different asymptotic densities (see also in
Section \ref{section:estimates}).
 We present below four generators of random normal forms.

\subsection{Generator of reduced forms}
The procedure $RG_{rf}$ generates a random element in a reduced
form of a syllable length $k$. This procedure depends on a given
probability distribution $\theta :\mathbb{N} \rightarrow
\mathbb{R}^+$ on the set of natural numbers $\mathbb{N}$ with
zero, two fixed probability distributions $\mu_A$ and $\mu_B$ on
$A\smallsetminus C$ and $B \smallsetminus C$, and two probability
distributions $\mu_{A,C}, \mu_{B,C}$ on $C$, where $C$ is viewed
as a subgroup of $A$ or $B$ correspondingly.

\begin{procedure}{(\em Generator {\bf $RG_{rf}$} of a random element in the reduced form (\ref{red_nf_without}))}\label{rgrf}\\
{\sc Input:} Number $k$ chosen with respect to a fixed
probability distribution $\theta :\mathbb{N} \rightarrow \mathbb{R}^+$.\\
{\sc Output:} A random word $u$ in the reduced form of length $k$.\\
{\sc Computations:}
\begin{itemize}

  \item [1)] Choose $A$ or $B$ with equal probability $\frac{1}{2}$.
    \item [2)]  If $k = 0$ then
\begin{itemize}
    \item [a)] if the choice in 1) is $A$ then choose randomly
  an element $c$ in $C$ with probability $\mu_{A,C}$;
   \item [b)] if the
  choice in 1) is $B$ then choose randomly
  an element $c$ in $C$ with probability $\mu_{B,C}$.
   \end{itemize}
 Output  $u = c.$
   \item [3)] If $k > 0$ then do the following
   \begin{itemize}
    \item [a)] if  the choice in 1) is $A$ then choose $g_1
   \in A\smallsetminus C$ with probability $\mu_A,$ then an element
   $g_2 \in B\smallsetminus C$ with probability $\mu_B$, and repeat
   this process choosing alternatively  $g_i \in A\smallsetminus C$ and $ g_{i+1} \in B \smallsetminus C$
     until $k$ elements
   $g_1, \ldots, g_k$ are constructed.
   \item [b)] if  the choice in 1) is $B$ then choose $g_1
   \in B\smallsetminus C$ with probability $\mu_B,$  then an element
   $g_2 \in A\smallsetminus C$ with probability $\mu_A$, and repeat this process as in step 3.a).
\end{itemize}
Output $u = g_1 \ldots g_k$.
\end{itemize}
\end{procedure}
\begin{remark} Generator {\bf $RG_{ef}$} of a random element $g$ in the freely
reduced form (\ref{any_g}) can be constructed in the similar way.
Namely, we should take elements from $A$ and $B$ consequently to
get a result. We will use this generator later in Section
\ref{section:estimates}.
\end{remark}

\subsection{Generator of canonical normal  forms}

Let $G=\mathop{A\ast B}\limits_C$, and  suppose $S, \, T$ are
fixed Schreier transversals for $C$ in $A$ and $B$ respectively.

Denote by $RG_{cnf}$  the following procedure for generating of
random elements in the canonical normal form of syllable length
$k.$ This procedure depends on a given probability distribution
$\theta :\mathbb{N} \rightarrow \mathbb{R}^+,$  two fixed
probability distributions $\mu_A$ and $\mu_B$ on $A\smallsetminus
C$ and $B \smallsetminus C$, and two probability distributions
$\mu_{A,C}, \mu_{B,C}$ on $C$.

\begin{procedure}\label{rgnf}{\em (Generator {\bf $RG_{cnf}$} of a random element in the canonical normal form (\ref{with}))}\\
{\sc Input:} A natural number $k$ chosen with respect to a fixed
probability distribution $\theta :\mathbb{N} \rightarrow \mathbb{R}^+$.\\
{\sc Output:} A random word $v$ in the canonical normal form of length $k$.\\
{\sc Computations:}
\begin{itemize}

  \item [1)] Choose $A$ or $B$ with equal probability
  $\frac{1}{2}$ and do as in the Procedure \ref{rgrf}:
    \begin{itemize}
    \item [a)] if the choice in 1) is $A$ then choose randomly
  an element $c$ in $C$ with probability $\mu_{A,C}$;
   \item [b)] if the
  choice in 1) is $B$ then choose randomly
  an element $c$ in $C$ with probability $\mu_{B,C}$;
   \end{itemize}
    \item [2)]  If $k = 0$ then output  $v = c$.

  \item [3)]  If $k \geq 1$ and

   \begin{itemize}
    \item [a)] the choice in 1) is $A$ then choose $g_1
   \in A \smallsetminus C ,$ represent it as $g_1 = c_1 s_1,$ where $c_1 \in C, \, s_1 \in S$ (so,
   $\mu_A(Cg_1)=\mu_A(Cs_1)$) and repeat this
   choosing alternatively  $g_i \in A\smallsetminus C$ and $ g_{i+1} \in B \smallsetminus C$
     with probabilities $\mu_A(Cg_i), \,\,\, \mu_B(Cg_{i+1})$ and represent
     $g_i = c_i s_i,\,\,\, g_{i+1} =c_{i+1}t_{i+1}$ until $k$ elements
   $s_1, t_2, s_3, t_4, \ldots $ are constructed.
\par Output $v = c s_1 t_2 s_3 t_4 \ldots $.

   \item [b)] the choice in 1) is $B$ then choose $g_1
   \in B \smallsetminus C$, represent it in a form $g_1 = c_1 t_1,$ where $c_1 \in C, \, t_1 \in T$  and repeat this procedure as in
   3a).

\par Output $v = c t_1 s_2 t_3 s_4 \ldots $.
\end{itemize}
\end{itemize}
\end{procedure}
\begin{remark}
Two probability distributions $\mu_A$ and $\mu_B$ on
$A\smallsetminus C$ and $B \smallsetminus C$ we describe in
details in Section \ref{subsection:measuringnfatomic}.
\end{remark}

Now we construct one more generator.

\subsection{Generator of cyclically reduced normal forms}

Let $RG_{crf}$ be the following procedure for generating of random
elements in the cyclically reduced canonical forms of length $2k$
or $1.$ This procedure depends on a given probability distribution
$\theta : 2\mathbb{N}\cup \{1\} \rightarrow \mathbb{R}^+,$ two
fixed probability distributions $\mu_A$ and $\mu_B$ on
$A\smallsetminus C_A^{\ast}$ and $B \smallsetminus C_B^{\ast}$,
and two probability distributions $\mu_{A,C}, \mu_{B,C}$ on $C,$
where $C_A^{\ast} = \mathop{\bigcup}\limits_{x\in A}C^{x}$ and
$C_B^{\ast} = \mathop{\bigcup}\limits_{x\in B}C^{x}.$

\begin{procedure}\label{rgcrf}
{\em (Generator {\bf $RG_{crf}$} of a random element in the cyclically reduced normal form (\ref{crf}))}\\
{\sc Input:} An even natural number $k$ or $1$ chosen with respect
to $\theta$.\\
{\sc Output:} A random word $w$ in the cyclically reduced normal form of length $k$.\\
{\sc Computations:}
\begin{itemize}

  \item [1)] Choose $A$ or $B$ with equal probability $\frac{1}{2}$ and do as in Procedure \ref{rgrf}:
  \begin{itemize}
    \item [a)] if the choice in 1) is $A$ then choose randomly
  an element $c$ in $C$ with respect to probability $\mu_{A,C}$;
   \item [b)] if the
  choice in 1) is $B$ then choose randomly
  an element $c$ in $C$ with respect to probability $\mu_{B,C}$;
   \end{itemize}
\item [2)]  If $k = 0$ then output $w = c.$
  \item [3)]  If $k = 1$
   \begin{itemize}
    \item [a)] if the choice in 1) is $A$ then choose randomly
  an element $g_1 \in A\setminus C_A^{\ast}$ with the probability $\mu_{A\setminus C_A^{\ast}}(Cg_1),$
  represent it as $g_1 = c_1 s_1,$ where $c_1 \in C, \, s_1 \in S$ and
    output $w= c s_1.$
   \item [b)] if the
  choice in 1) is $B$ then choose randomly
  an element $g_1 \in B\setminus C_B^{\ast}$ with the probability $\mu_{B\setminus
  C_B^{\ast}}(Cg_1),$  represent it as $g_1 = c_1 t_1,$ where $c \in C, \, t_1 \in T$ and
output $w= c t_1.$
   \end{itemize}
   \item [4)] If $k = 2l, \,\,\,l \geq 1$ then do the following
   \begin{itemize}
    \item [a)] if the choice in 1) is $A$ then choose $g_1
   \in A\smallsetminus C^{\ast}$ as in 3), then an element $g_2
   \in B\smallsetminus C^{\ast}$ with probability $\mu_B(Cg_2)$ and represent $g_2 = c_2t_2,$ and repeat
   this process choosing alternatively  $g_i \in A \smallsetminus C^{\ast}$ and $g_{i+1} \in B \smallsetminus C^{\ast}$
     with probabilities $\mu_A(Cg_i), \,\,\, \mu_B(Cg_{i+1})$ and represent
     $g_i = c_i s_i,\,\,\, g_{i+1} =c_{i+1}p_{i+1}$ until $k$ elements
   $s_1, t_2, s_3, t_4, \ldots $ are constructed.

Output $w = c s_1 t_2 s_3 t_4 \ldots $.

   \item [b)] if  the choice in 1) is $B$ then choose $g_1
   \in B \smallsetminus C^{\ast}$ as in 3), then an element $g_2
   \in A \smallsetminus C^{\ast}$ with probability $\mu_A(Cg_2)$ and represent $g_2 = c_2s_2,$ and repeat as in
   4a).

Output $w = c t_1 s_2 t_3 s_4 \ldots $.

\end{itemize}
\end{itemize}
\end{procedure}
\begin{remark}
We shall talk about probability distributions on cosets
$\mu_A(Cg_i)$, $\mu_B(Cg_i)$ and other sets later in Section
\ref{section:estimates}.
\end{remark}

\section{Evaluation of sets of randomly generated normal forms}\label{section:estimates}

Let $G  = \mathop{A\ast B}\limits_{C}$, where  $A$ and $B$ are
free groups with finite bases $X$ and $Y$ correspondingly, and
 $C$ is a finitely generated subgroup. Now we are ready to stratify sets
of generated random normal forms of elements in $G$ into regular,
singular and stable or unstable subsets and evaluate their sizes
in whole sets of corresponding forms.
\subsection{Atomic measures of elements in normal forms}\label{subsection:measuringnfatomic}
In this section we introduce probability measures on sets of
normal forms, constructed with a help of generators above.

The probability to obtain an element $g \in G$ in freely reduced
or reduced form of syllable length $k$ on the output of generator
$RG_{ef}$ or $RG_{rf}$ is equal to
$$\mu_k(g) = \frac{1}{2} \mu_1(g_1)\ldots\mu_k(g_k),$$
 where $\mu_i(g_i) = \mu_A(g_i)$ if $g_i \in A$, and $\mu_i(g_i) =
 \mu_B(g_i)$ otherwise, $i = 1, \ldots, k;$ $k \geq 1;$
 and if $k=0$ then $\mu_0(c) = \mu_C(c),$ where $\mu_C(c) =
 \mu_{A,C}$ if $C$ is viewed as a subgroup of $A$, and $\mu_C(c) =
 \mu_{B,C}$ otherwise. Clearly, $\mu_k$ is an atomic probability
 measure on sets ${\mathcal EF}_k$ or ${\mathcal RF}_k$ of all (freely) reduced elements of length $k.$
 Now one can calculate a probability measure
 $\mu$ on sets ${\mathcal EF}, {\mathcal RF}$  of all (freely) reduced
 elements:
  $$\mu(g) = \theta(k)\mu_k(g).$$

To complete definitions above we have to define probability
measures $\mu_A$, $\mu_B$, $\mu_{A,C}$, and $\mu_{B,C}$. There are
a lot of different methods to describe a probability distribution
on a free group; for example, it can be done as in \cite{fmrI}
with a help of no-return random walk $W_s$ ($s \in (0,1]$) on the
Cayley graph of $A = F(X)$ of rank $r = |X|$. The probability
$\mu_s(g)$ for this process to terminate at $g$ is given by the
formula
\begin{equation}
\mu_s(g) = \frac{s(1-s)^{|g|}}{2m\cdot (2m-1)^{|g|-1}} \quad
\hbox{ for } w \ne 1
\end{equation}
and
\begin{equation}
\mu_s(1) =s.
\end{equation}
This random walk can be considered on the set $A \smallsetminus C$
with small changes. Since the set $A \smallsetminus C$ is regular
in $A = F(X)$ there is a simple procedure to define probability
measures using random walks in the corresponding finite automata
(or graph). We introduce these measures using different way, i.e.
infinite trees, the Cayley graphs of $A$ and $B$. Let $\Gamma$ be
the Cayley graph of $A$ with respect to a basis $X$ of $A,$
$\Gamma_C$ is a subgroup graph for $C$, and $\Gamma_C^{\ast}$ its
extended graph. Denote by $\pi : \Gamma \rightarrow
\Gamma_C^{\ast}$ the unique canonical projection from $\Gamma$
onto $\Gamma_C^{\ast}$, so $\pi$ is a morphism of graphs which
preserves labels. We choose a real number $s \in (0,1)$ and define
a no-return random walk ${\mathcal W}_s$ on $\Gamma$ as follows.
If ${\mathcal W}_s$ is at vertex $v \in \Gamma$ and $v \not \in V
(\Gamma_C)$ then ${\mathcal W}_s$ stops at $v$ with probability
$s$ and moves from $v$ (away from the root $1_C$) along an
adjacent edge with equal probability $\frac{1-s}{2r-1}$.  If
${\mathcal W}_s$ is at vertex $v \in V (\Gamma_C)$  then
${\mathcal W}_s$ moves from $v$ (away from the root $1_C$) along
any adjacent edge with equal probability $\frac{1}{2r-1}$. The
random walk ${\mathcal W}_s$ induces (via the projection $\pi$) a
random walk ${\mathcal W}_s^{\ast}$ on $\Gamma_C^{\ast}$. Then the
probability $\mu_s(w)$ for ${\mathcal W}_s$ to stop at $w \in A
\setminus C$ can be written  as follows
 $$\mu_s(w) = \frac{1}{2r}\frac{(1-s)^{\abs{w}-m_w}}{(2r-1)^{\abs{w}-m_w}}
 \frac{1}{(2r-1)^{m_w -1}},$$
 where $\abs{w}$ is a freely-reduced length of $w$ in $A$
 and $m_w$ is the number of times the path $\pi(w)$ visits the
 vertex $1_C$ in $\Gamma_C^{\ast}.$

Now we calculate probability to obtain an element in the canonical
and cyclically reduced normal forms. Namely, the probability to
obtain $v = c p_1 p_2 \ldots p_k,$ $k \geq 0$, is equal to
$$\mu_k(v)=\frac{1}{2}\mu_C(c) \mu_1(p_1)\mu_2(p_2)\ldots \mu_k(p_k),$$
where $\mu_i(p_i)= \mu_{F(p_i)}(Cp_i),$ $i = 0, \ldots, k$ and
$\mu_C(c) =  \mu_{A,C}$ if $C$ is viewed as a subgroup of $A$, and
$\mu_C(c) = \mu_{B,C}$ otherwise.

Clearly, $\mu_k$ is an atomic probability
 measure on the set $\mathcal{CNF}_k$ of all canonical normal (or $\mathcal{CRF}_k$  of cyclically reduced) forms of syllable length $k.$
 Then the probability measure
 $\mu$ on sets $\mathcal{CNF}, \mathcal{CRF}$ of all canonical (or cyclically reduced) normal forms
 is equal to
  $$\mu(v) = \theta(k)\mu_k(v).$$

To complete this definition, we have to know how to calculate
probability measures of subgroup $C$ and its cosets. It can be
done with a help of consolidated subgroup graph for $C$ (see
Section 3.3 of \cite{multiplicative} for details).

\subsection{Measures on sets of normal forms}\label{subsection:reg_stable}
In the rest of the paper by $\mathcal{NF}_r (\mathcal{NF}_s)$ we
denote the set of all regular (stable) elements in normal form
(prefix $\mathcal{NF}$ will changes depends on type of chosen
form) and by $\mathcal{NF}_{\rm sin} (\mathcal{NF}_{\rm uns})$ its
complement in $\mathcal{NF}.$

Decision of algorithmic problems and problem of stratification of
inputs leads to necessity of estimation of sizes of normal forms
and their regular and stable subsets. For this purpose we use  two
following approaches: asymptotic approach and related notion of
$L-$measure and probabilistic approach and the notion of
$\lambda_L-$measure. To apply classification theorems about
regular sets, we describe in this section structure of
$\mathcal{NF}$ and $\mathcal{NF}_{\rm uns}$.

Remind some definitions first. In \cite{fmrI} for subsets $R, L$
of a free group $F$ of a finite rank was defined their {\em size
ratio} at length $n$ by
   $$f_n(R,L) = \frac{f_n(R)}{f_n(L)} = \frac{|R\cap S_n|}{|L\cap S_n|}.$$
  The {\em asymptotic density} of $R$ relative to $L$ is defined by
  $$\rho(R,L) = \overline{\lim\limits_{n \rightarrow \infty}} f_n(R,L).$$
  By $r_L(R)$ we denote the {\em cumulative size ratio} of $R$ relative to $L$:
  $$  r_L(R) = \sum_{n=1}^\infty f_n(R,L).  $$
  $R$ is called {\em negligible} relative to $L$  if $\rho(R,L) = 0$.
   A set $R$ is termed {\em exponentially negligible} relative to  $L$ (or  {\em exponentially $L$-negligible}) if
   $f_n(R,L) \leq \delta^n$ for all sufficiently large $n$ and some positive constant $\delta <1$.
   Further, a set $R$ is called {\em generic} relative to $L$  if
   $\lim\limits_{n \rightarrow \infty} f_n(R,L)$ exists and equal to $1$; $R$ is termed {\em exponentially generic}
   relative to $L$  if there exists a positive constant $\delta <1$
such that $1-\delta^n < f_n(R,L) < 1$ for large enough $n$. We
will use in this section the following notion of a thick set;
namely, a set $R$ is {\em thick} relative to $L$ if
$\lim\limits_{n \rightarrow \infty} f_n(R,L)$ exists and strictly
greater than $0$.

In \cite{fmrI} were also introduced notions of
$\lambda_L-$measurable and exponentially $\lambda_L-$measurable
set $R$. This way of measuring closed to the frequency measure and
also coming up from no-return non-stop random walk on an automaton
for $L$. Namely,
$$\lambda_L(R)  =  \sum_{w \in R}\lambda_L(w) = \sum_{n=0}^{\infty}f^{\prime}_n(R,L),$$
where $$f^{\prime}_n(R,L) = \sum_{w \in R \cap S_n}
\lambda_L(w),$$ and $\lambda_L(w), f^{\prime}_n(R,L)$ defined by a
random walk on an automaton for $L$ (see Section 5.2, p. 109 of
\cite{fmrI} for details). Note, that generally speaking,
$f^{\prime}_n(R,L)$ differs from $f_n(R,L)$ defined above because
of specific of an automaton for $L$. On the other hand,
$f^{\prime}_n(R,F(X)) = f_n(R,F(X))$. Notice also, that
$\lambda_L$ is multiplicative, since $\lambda_L(u v) =
\lambda_L(u) \lambda_L(v)$ for any $u,v \in R$ such that $uv = u
\circ v$ and $uv \in R$; it is easy to check by definition of
random walk on $L$.

  We say that $R$ is {\em $\lambda_L$-measurable}, if $\lambda_L(R)$ is finite.
   A set $R$ is termed {\em exponentially $\lambda_L$-measurable}, if
   $f_n^{\prime}(R,L) \leq q^n$ for all sufficiently large $n$ and some positive constant $\delta <1$.

In \cite{fmrI} was given the following definition. For every $w
\in F$ the set $C_L(w) = L \cap C(w)$ is called an $L$-cone, and
$C_L(w)$ is called \emph{$L$-small}, if it is exponentially
$\lambda_L$-measurable.

\begin{theorem}\label{Lcones}\cite[Theorem 5.4]{fmrI}.
Let $R$ be a regular subset of a prefix-closed regular set $L$ in
a finite rank free group $F$. Then either the  prefix closure
$\overline{R}$ of $R$ in $L$ contains a non-small  $L-$cone or
$\overline{R}$ is exponentially $\lambda_L$-measurable.
\end{theorem}
This theorem is very convenient for asymptotic estimates of
subsets of $F$ and we are going to apply it in Theorem
\ref{nf_lambdaLmeasurable}. To prove this theorem  we will use
also lemmata \ref{le:f_0expnegl} and \ref{le:fp_Lcones}:

\begin{lemma}\label{le:f_0expnegl} Suppose $F=F(X)$ is a finite rank free group and $W = \{w_1, \ldots, w_d
\}$ is a finite subset of $F$ such that at least one $w_i$ is not
trivial. Then
$$F_0 = \{ f \in F| f \textrm{ does not contain any element of $W$
as a subword } \}$$ is a regular, prefix-closed exponentially
$\lambda_F-$measurable set.
\end{lemma}

\proof Remind Myhill-Nerode criterion, which gives necessary and
sufficient conditions on a subset of semigroup to be regular. For
a language $R$ in semigroup $A^{\ast}$ consider an equivalence
relation $\mathop{\sim}\limits_{R}$ in $A^{\ast}$ relative to $R$:
$v_1 \mathop{\sim}\limits_{R} v_2$ if and only if for each string
$u$ over $A$ the words $v_1 u$ and $v_2 u$ are either
simultaneously  in $R$ or not in $R$.

\begin{quote}{\bf{Myhill-Nerode.}} A set $R$ is regular in
$A^{\ast}$ iff there are only finitely many
$\mathop{\sim}\limits_{R}-$equivalence classes.
\end{quote}

Proof of this theorem can be found, for example, in \cite{eps}
(see Theorem 1.2.9).

Since we works in a free group, the analogue of latter theorem for
groups was proved in \cite{fmrI}. Now, let  $R \subseteq F$.
Define an  equivalence relation $\sim$ on $F$ relative to $R$ such
that $v_1 \mathop{\sim}\limits_{R} v_2$ if and only if for each $u
\in F$ the following condition holds: $v_1u = v_1 \circ u$ and
$v_1u \in R$ if and only if $v_2 u = v_2 \circ u$ and $v_2 u \in
R.$

\begin{lemma}\label{Myh-Ner_group}\cite[Lemma 5.3]{fmrI}.
Let  $R \subseteq F$. Then $R$ is regular if and only if there are
only finitely many $\mathop{\sim}\limits_{R}-$equivalence classes.
\end{lemma}

It is sufficient to show the statement of lemma
\ref{le:f_0expnegl} for a singleton set $W$. Indeed, let $F_0$ be
the subset of words in $F$ that do not contain $w$ as a subword,
and $F_1$ be the set of words that do not contain $w, w_1, \ldots,
w_d$. Then $F_1 \subseteq F_0$, and $F_1$ is exponentially
$\lambda_F-$measurable if $F_0$ is so. Assume now that $W = \{ w
\}$ and $|w| = t \ge 1$.

We describe all equivalence classes relative to $F_0$:

\bi \item [1)] $\mathcal{K}_0=  F \smallsetminus F_0$;

\item [2)] for all $u \in S_{t-1}$ define $\mathcal{K}_{u}= \{ f
\in F: f = g \circ u, g \in F, u \in S_{t-1} \} \smallsetminus
\mathcal{K}_0$;

\item [3)] for all $v \in B_{t-2}$ set $\mathcal{K}_{v}= \{ v \}$.

\ei

Now we will prove that this decomposition define an equivalence
relation on $F$ relative to $F_0$, i.e. for every representative
$u_1, u_2 \in F$ holds $u_1 \mathop{\sim}\limits_{F_0} u_2$
$\Leftrightarrow$ (for all $p \in F:$ $u_1 p = u_1 \circ p$ and
$u_1p \in F_0$ iff $u_2 p = u_2 \circ p$ and $u_2 p \in F_0$).

Obviously, this condition holds for classes $\mathcal{K}_0$ and
$\mathcal{K}_{v}$ from i.3. We shall show that it holds also for
classes $\mathcal{K}_{v}$. Let $u_1, u_2 \in \mathcal{K}_{u}$,
where $u \in S_{t-1}$, i.e. elements $u_1 = g_1 \circ u, u_2 = g_2
\circ u$ and $u_1, u_2$ do not contain $w$ as a subwords. If $u_1
p = (g_1 \circ u) \circ p \in F_0$, then $u_2 p = g_2 \circ u
\circ p$ by definition of class. So we should show also, that $g_2
\circ u \circ p$ doesn't belong to $\mathcal{K}_0$. But since
neither $g_2 \circ u$ nor $u \circ p$ contains $w$ as a subwors,
then $w$ should have a form $w = g_2' \circ u \circ p'$, where
nontrivial element $g'_2$ is an end of $g_2$ and $p'$ is
nontrivial beginning  of $p$. But it implies that $t = |w| =
|g_2'| + |u| + |p'|
> t$, a contradiction.
Suppose now that $u_1 p = (g_1 \circ u) \circ p \notin F_0$. Then
exactly $u \circ p$ contains $w$ as a subword. Hence, the element
$u_2 p = g_2 \circ u \circ p$ also contains $w$ as a subword.

Therefore, $\mathop{\sim}\limits_{F_0}$ is a relation equivalence
on $F$ relative to $F_0$,and since the number of equivalence
classes relative to $F_0$ is finite, it follows from lemma
\ref{Myh-Ner_group} that $F_0$ is regular in $F$. By definition of
$F_0$ it is prefix-closed in $F$. Further, since $F_0$ evidently
doesn't contain a cone, then by Theorem \ref{Lcones} $F_0$ is
exponentially $\lambda_F-$measurable. The last statement about
exponentially measurability of $F_0$ in $F$ follows also from
Lemma 3 in \cite{AC}. \eproof

\begin{lemma}\label{le:fp_Lcones} Suppose $F=F(X)$ is a finite rank free group and $W = \{w_1,
\ldots, w_d \}$ is a finite subset of $F$, at least one of $w_i$
is not trivial and $F_0$ is a set of words that do not contain
elements of $W$ as subwords. Let $F_1$ denote the supplement of
$F_0$ in $F$. Let $R \subseteq F_0$ be regular and $L \subseteq F$
be regular prefix-closed in $F$. If for every non-small cone
$C_L(u)$ in $L$ holds $C_L(u) \cap F_1 \neq \emptyset$, then $R$
is exponentially $\lambda_L-$measurable.
\end{lemma}
\proof By assumption of lemma the set $R \cap F_1$ is nonempty,
but all non-small cones $C_L(u)$ have nonempty intersetion with
$F_1$; therefore, $C_L(u) \nsubseteq R$ and by theorem
\ref{Lcones} the set $R$ is exponentially
$\lambda_L-$measurable.\eproof

To estimate sizes of subsets of normal forms we prove first the
following theorem which describes their structure.

\begin{theorem}\label{nf_structure}
Let $G=\mathop{A\ast B}\limits_C$ be an amalgamated product, where
$A,B,C$ are free groups of finite rank. Then the set
$\mathcal{NF}_{\rm uns}$ is regular in $G$ and the set
$\mathcal{NF}$ is regular prefix closed in $G$  for all
$\mathcal{NF} = \{ \mathcal{EF}, \mathcal{RF}, \mathcal{CNF},
\mathcal{CRF} \}$.
\end{theorem}

\proof Suppose $L,M$ are two regular sets in alphabets $X\cup
X^{-1},Y \cup Y^{-1}$ correspondingly (recall that $A=F(X), B=
F(Y)$). Denote by $LM$ concatenation of sets $L$ and $M$ and by
$L^{\ast}$ a monoid generated by $L$. We will often use the
following formula in the sequel:
\begin{equation}\label{free_product}
L \ast M = L \sqcup M \sqcup LML \sqcup MLM \sqcup L(ML)^{\ast}M
\sqcup M(LM)^{\ast}L
\end{equation}
Particulary, it shows that for regular $L,M$ the set $L\ast M$ is
also regular (in alphabet $X\cup X^{-1}\cup Y \cup Y^{-1}$).

Consider the set of all freely reduced normal forms:

\begin{equation}\label{nf_ef}
\mathcal{EF} = A \ast B.
\end{equation}
Evidently, the set $\mathcal{EF}$ is regular and prefix closed.

The set of all unstable freely reduced normal forms:
\begin{equation}\label{nf_ef_uns}
\mathcal{EF}_{\rm uns} = (\mathop{\cup}\limits_{s \in S_{\rm
uns}}Cs) \ast (\mathop{\cup}\limits_{t \in T_{\rm uns}}Ct).
\end{equation}

By Lemma \ref{pr:basic-properties} the set $S_{\rm uns}$ is a
finite union of  left cosets of $C$ of the type $s_1s_2^{-1}C$,
where $s_1, s_2 \in S_{\rm int}$. The set $\mathcal{EF}_{\rm uns}$
is regular as a free product of regular sets (see formula
(\ref{free_product})) since $\mathop{\cup}\limits_{s \in S_{\rm
uns}}Cs$ (as well as $\mathop{\cup}\limits_{t \in T_{\rm uns}}Ct$)
is a concatenation of regular sets $C$ and $S_{\rm uns}$ (or
$T_{\rm uns}$).

Consider all non-trivial reduced forms:
\begin{equation}\label{nf_rf}
\mathcal{RF} = (A \setminus C) \ast (B \setminus C)
\end{equation}

The set of all unstable reduced forms:
\begin{equation}\label{nf_rf_uns}
\mathcal{RF}_{\rm uns} = (A_1' \setminus C) \ast (B_1' \setminus
C)
\end{equation}
such that $$A_1' = \mathop{\cup}\limits_{s \in S_{\rm uns}}Cs
\,\,\textrm{ and } \,\,B_1' = \mathop{\cup}\limits_{t \in T_{\rm
uns}}Ct.$$

Since difference of regular sets is regular again, both sets of
forms are evidently regular. To see that $A \setminus C$ is a
prefix closed in $A$ set, one can identify it with a language of
all words in a Schreier graph $\Gamma^{\ast} = \Gamma_C^{\ast}$
readable as a labels of paths starting in a root vertex $1_C$ (and
probably return to this vertex again) but finish in arbitrary
vertex of $\Gamma^{\ast}$ except the root one (remind that the
graph $\Gamma^{\ast}$ was defined in Section
\ref{subsection:measuringS}). Clearly, such a language is prefix
closed in $A$. Analogously, $B \setminus C$ is prefix closed in
$B$ and so $\mathcal{RF}$ is prefix closed as a free product of
prefix closed sets.

 Let $\mathcal{NF}$
be the set of all canonical normal forms. Every $v \in
\mathcal{CNF}$ can be written in the form (\ref{with}), i.e. as
$$v = c p_1 p_2 \ldots p_l,$$ where $c \in C, \,\, p_i \in
(S \cup T) \setminus 1, \textrm{ and } F(p_i) \neq F(p_{i+1}), \,
i = 1, \ldots l, l \geq 0.$ Now consider the set of all canonical
normal forms:

\begin{equation}\label{nf_cnf}
\mathcal{CNF} = C \mathop{\circ}\limits_t ( S \ast T)
\end{equation}

At the same time,

\begin{equation}\label{nf_cnf2}
\mathcal{CNF} = A \circ (S \ast T)_1 \sqcup B \circ (S \ast
T)_2,
\end{equation}
 where $(S \ast T)_1$ starts from $t \in T$ and $(S \ast
T)_2$ starts from $s \in S.$ Substitute $\circ$ with concatenation
and applying formula (\ref{free_product}) for $(S \ast T)_i,
i=1,2$, conclude that $\mathcal{CNF}$ is regular set. Since $A,S$,
and $T$ are prefixed closed, the set $\mathcal{CNF}$ also has this
property.

Analogous decomposition can be applied to $\mathcal{CNF}_{\rm
uns}:$
\begin{equation}\label{nf_cnf_uns}
\mathcal{CNF}_{\rm uns} =  C \mathop{\circ}\limits_t ( S_{\rm uns}
\ast T_{\rm uns}) = A_1' \circ (S_{\rm uns} \ast T_{\rm uns})_1
\sqcup B_1' \circ (S_{\rm uns} \ast T_{\rm uns})_2,
\end{equation}

where $A'_1, B'_1$ obtained as intersections of two regular sets
and $(S_{\rm uns})_i, (T_{\rm uns})_i, i =1 , 2$ are regular by
Proposition \ref{pr:basic-properties}.

The proof of this theorem for $\mathcal{CRF}$ is straightforward.
\eproof

\subsection{Regular and stable normal forms}\label{subsection:reg_stable}

Notions of regular and stable forms was formulated in Section
\ref{subsection:measuringS}, and now we are ready to prove the
main theorems about evaluation of the set of singular and unstable
(and, therefore, regular and stable) normal forms in group $G.$
\begin{maintheorem}\label{nf:basic}
Let $G=\mathop{A\ast B}\limits_C$ be an amalgamated product, where
$A,B,C$ are free groups of finite rank. Then for every set of
normal forms $\mathcal{NF} = \{ \mathcal{EF}, \mathcal{RF},
\mathcal{CNF}, \mathcal{CRF} \}$

\bi

\item [(i)] If $C$ has a finite index in $A$ and in $B,$ then
every normal form is singular and unstable, i.e. $\mathcal{NF}_{\rm
sin} = \mathcal{NF}_{\rm uns} = \mathcal{NF};$

\item [(ii)] If $C$ of infinite index either in $A$ or in $B,$
then $\mathcal{NF}_r$ and $\mathcal{NF}_s$ are exponentially
$\mu-$generic relative to $\mathcal{NF},$ and $\mathcal{NF}_{\rm
sin}$ and $\mathcal{NF}_{\rm uns}$ are exponentially
$\mu-$negligible relative to $\mathcal{NF}$ in the following cases:

\bi \item [(ii.1)] $\mu$ is defined by pseudo-measures $\mu_A$ and
$\mu_B$, which are cardinality functions on $A$ and  $B$
correspondingly; in this case $\rho_{\mu}$ is a bidimensional
asymptotic density;

\item [(ii.2)] $\mu$ is defined by atomic probability measures
$\mu_{A,l}$ and $\mu_{B,l}$ on $A$ and $B$ correspondingly; in
this case $\rho^C$ is a bidimensional Cesaro asymptotic density.
\ei \ei

\end{maintheorem}

\proof Observe, that all singular representatives by Proposition
\ref{pr:basic-properties} are also unstable, and it will be
sufficient to prove Theorem \ref{nf:basic} for unstable normal
forms only.

Suppose first that $C$ has a finite index in both $A$ and $B$.
Then there are nontrivial subgroups $N_A$ in $C$, which is normal
subgroup of $A$, and $N_B$ in $C$, which is normal in $B$.
Therefore, $A \ast B = N^*_{A \ast B}(C).$ Then by definition all
cosets representatives (and hence $\mathcal{NF}$) are unstable,
i.e. $\mathcal{NF}_{\rm sin} = \mathcal{NF}_{\rm uns} =
\mathcal{NF}.$

Suppose now, that subgroup $C$ has an infinite index in $A$. Let
$\mathcal{NF} = \{ \mathcal{EF}, \mathcal{RF}\}$ and (ii.1) is
hold. By Proposition \ref{munegligibility} it is sufficient to
show that $\mathcal{NF} = A_0 \ast B_0, \,\, \mathcal{NF}_{\rm
uns} = A_1 \ast B_1$ and $A_1$ is exponentially negligible
relative to $A_0$ (it is clear that the density
$\rho_{\mu_B}(B_1,B_0)$ for these forms exists because of
definitions of $B_0$ and $B_1$).

Consider freely reduced forms first. Applying formulae
(\ref{nf_ef}), (\ref{nf_ef_uns}), set $A_0 = A; B_0= B$ and $A_1 =
\mathop{\cup}\limits_{s \in S_{\rm uns}}Cs; B_1 =
\mathop{\bigcup}\limits_{t \in T_{\rm uns}}Ct.$ Observe, that
$S_{\rm uns}$ is exponentially negligible in $S$ by Corollary
\ref{s_nst-ins} and so $A_1$ is exponentially negligible in $A_0$
by Proposition 4.7. from \cite{fmrI}.

For the set of all non-trivial reduced forms, using (\ref{nf_rf}),
(\ref{nf_rf_uns}), set $A_0 = A \setminus C$; $B_0= B \setminus C$
and $A_1 = A_1' \setminus C$; $B_1 = B_1' \setminus C$, where
$A_1' = \mathop{\cup}\limits_{s \in S_{\rm uns}}Cs \,\,\textrm{
and } \,\,B_1' = \mathop{\cup}\limits_{t \in T_{\rm uns}}Ct.$ We
have already shown for freely reduced forms, that $A_1'$ is
exponentially negligible in $A,$ and hence there is a $q, \,\,0 <
q <1,$ such that $\fracd{\mu_A((A_1')_n)}{\mu_A((A)_n)} < q$ for
all $n \geq n_0$ for some natural number $n_0.$

It is clear that $\fracd{\mu_A((C)_n)}{\mu_A((A)_n)} < q,$ where
$n \geq n_0.$ Then $$\mu_A((A_1)_n)= \mu_A(((A_1')_n) \setminus
((C)_n))=$$ $$=\mu_A((A_1')_n) - \mu_A((C)_n) < q \mu_A((A)_n) -
\mu_A((C)_n)<$$
$$<q (\mu_A((A)_n) - \mu_A((C)_n)).$$ Then we obtain
$\fracd{\mu_A((A_1)_n)}{\mu_A((A_0)_n)} < q$ for all $n \geq n_0$
and so the set $A_1$ is exponentially negligible relative to
$A_0.$

Suppose now that the case (ii.2) holds, i.e. the measure $\mu$ is
defined on $F = A \ast B$ by atomic probability measures
$\mu_{A,l}$ and $\mu_{B,l}$ on $A$ and on $B$ correspondingly and
$\rho^C$ is a bidimensional Cesaro asymptotic density. Since $A_1'
= \mathop{\cup}\limits_{s \in S_{\rm uns}}Cs$ is exponentially
negligible in $A,$ by Lemma \ref{snegl_is_Cnegl} there exists a
number $q, \textrm { where } 0 < q < 1$  such that
$\fracd{\mu_{A,l}(A_1')}{\mu_{A,l}(A)} < q$ for all $l \geq l_0$.
Then by Proposition \ref{Cnegligibility} the set of all unstable
freely reduced forms is $C-$negligible relative to the set
$\mathcal{EF}$.

Further, since $\mu_{A,l}(A_1') - \mu_{A,l}(C) < q \mu_{A,l}(A) -
\mu_{A,l}(C) < q ( \mu_{A,l}(A) - \mu_{A,l}(C))$ for all $l \geq
l_0$, obtain
$$ \fracd{\mu_{A,l} (A_1' \setminus C)}{\mu_{A,l} (A \setminus C)} < q. $$
Therefore, by Proposition \ref{Cnegligibility} we obtain claim
(ii.2) for $\mathcal{RF}$.

Now consider canonical normal forms, given by formulae
(\ref{nf_cnf2}) and (\ref{nf_cnf_uns}):
$$\mathcal{CNF} = (A \circ (S \ast T)_1) \sqcup (B \circ (S \ast T)_2) \textrm{  and  } \mathcal{CNF}_{\rm uns} =
(A \circ (S_{\rm uns} \ast T_{\rm uns})_1) \sqcup (B \circ (S_{\rm
uns} \ast T_{\rm uns})_2).$$ Clearly, it is enough to show the
statement of the theorem for a pair $\Sigma=(A \circ (S \ast
T)_1)$ and $\Sigma_{\rm uns}=(A \circ (S_{\rm uns} \ast T_{\rm
uns})_1)$. Denote by $\widehat{A}$ the set $A \setminus \{ 1 \}$.

We shall show that unstable canonical normal forms are
$\mu-$exponentially negligible in the set of canonical normal
forms relative to bidimensional asymptotic density defined by
cardinality functions. By definition of a frequency function
$$\rho_{\mu}^{n,k}(\Sigma_{\rm uns}, \Sigma) = \fracd{\mu
(\Sigma_{\rm uns} \cap \Sigma_{n,k})}{\mu(\Sigma_{n,k})}=$$

$$=\fracd{\mu_A((\widehat{A})_{\leq n}) \mu(((S_{\rm uns} \ast
T_{\rm uns})_1)_{n,{k-1}})+\mu(((S_{\rm uns} \ast T_{\rm
uns})_1)_{n,k})}{\mu_A((\widehat{A})_{\leq n}) \mu(((S \ast
T)_1)_{n,{k-1}})+\mu(((S \ast T)_1)_{n,k})}=$$

$$=\fracd{\mu_A((\widehat{A})_{\leq n}) \mu(((S_{\rm uns} \ast
T_{\rm uns})_1)_{n,{k-1}})+\mu(((S_{\rm uns} \ast T_{\rm
uns})_1)_{n,k})}{\mu_A((\widehat{A})_{\leq n}) \mu(((S \ast
T)_1)_{n,{k-1}})+\mu(((S \ast T)_1)_{n,k})}= $$

$$ = \fracd{\mu_A((\widehat{A})_{\leq n})\cdot (\mu_B((T_{\rm
uns})_{\leq{n}}))^{[\frac{k}{2}]}(\mu_A((S_{\rm
uns})_{\leq{n}}))^{[\frac{k-1}{2}]}+(\mu_B((T_{\rm
uns})_{\leq{n}}))^{[\frac{k+1}{2}]}(\mu_A((S_{\rm
uns})_{\leq{n}}))^{[\frac{k}{2}]}}{\mu_A((\widehat{A})_{\leq
n})\cdot
(\mu_B((T)_{\leq{n}}))^{[\frac{k}{2}]}(\mu_A((S)_{\leq{n}}))^{[\frac{k-1}{2}]}
+(\mu_B((T)_{\leq{n}}))^{[\frac{k+1}{2}]}(\mu_A((S)_{\leq{n}}))^{[\frac{k}{2}]}}.$$

Suppose $k$ is odd. Therefore, $\rho_{\mu}^{n,k}(\Sigma_{\rm
uns},\Sigma) =$
$$ \fracd{\left(\mu_A((\widehat{A})_{\leq n}) +
\mu_B((T_{\rm uns})_{\leq{n}}))\right)(\mu_B((T_{\rm
uns})_{\leq{n}}))^{\frac{k-1}{2}}(\mu_A((S_{\rm
uns})_{\leq{n}}))^{\frac{k-1}{2}}}{\left(\mu_A((\widehat{A})_{\leq
n}) +
\mu_B((T)_{\leq{n}})\right)(\mu_B((T)_{\leq{n}}))^{\frac{k-1}{2}}(\mu_A((S)_{\leq{n}}))^{\frac{k-1}{2}}}.$$

Since set $S_{\rm uns}$ is also exponentially negligible relative
to $S$, $\rho_{\mu_A}(S_{\rm uns},S) < 1.$ Then

$\mathop{\lim}\limits_{n \rightarrow \infty} \fracd{\mu_A((S_{\rm
uns}))_n}{\mu_A((S)_n)} = \mathop{\lim}\limits_{n \rightarrow
\infty} \fracd{\mu_A((S_{\rm uns}))_{\leq n}}{\mu_A((S)_{\leq
n})}$ and $\fracd{\mu_A((S_{\rm uns}))_{\leq n}}{\mu_A((S)_{\leq
n})} < (1-\varepsilon_1)$.

By the same reason $\fracd{\mu_B((T_{\rm uns}))_{\leq
n}}{\mu_B((T)_{\leq n})} < (1-\varepsilon_2)$ for some $0 <
\varepsilon_1, \varepsilon_2 \leq 1$ and large enough $n$.

Thereby, $$\rho_{\mu}^{n,k}(\Sigma_{\rm uns},\Sigma) <
\left(\fracd{\mu_A((\widehat{A})_{\leq n}) + \mu_B((T_{\rm
uns})_{\leq{n}})}{\mu_A((\widehat{A})_{\leq n}) +
\mu_B((T)_{\leq{n}})}\right)\left((1 - \varepsilon_1)(1-
\varepsilon_2)\right)^{\frac{k-1}{2}}$$ for large enough $n$ and
odd $k$.

One can check that for all even $k$ frequencies
$\rho_{\mu}^{n,k}(\Sigma_{\rm uns},\Sigma)$ are bounded by
$\left(\fracd{\mu_A((\widehat{A})_{\leq n}) + \mu_A((S_{\rm
uns})_{\leq{n}})}{\mu_A((\widehat{A})_{\leq n}) +
\mu_A((S)_{\leq{n}})}\right)(1- \varepsilon_2)\left((1 -
\varepsilon_1)(1- \varepsilon_2)\right)^{\frac{k-2}{2}}.$

Fractions $\left(\fracd{\mu_A((\widehat{A})_{\leq n}) +
\mu_B((T_{\rm uns})_{\leq{n}})}{\mu_A((\widehat{A})_{\leq n}) +
\mu_B((T)_{\leq{n}})}\right)$ and
$\left(\fracd{\mu_A((\widehat{A})_{\leq n}) + \mu_A((S_{\rm
uns})_{\leq{n}})}{\mu_A((\widehat{A})_{\leq n}) +
\mu_A((S)_{\leq{n}})}\right)$ do not depend on $k$ and less than
$1$ for all large enough $n$.

 Let $\varepsilon=
\min{\{\varepsilon_1,\varepsilon_2\}}$. Then
$$\rho_{\mu}^{n,k}(\Sigma_{\rm uns},\Sigma) <
\left(\fracd{\mu_A((\widehat{A})_{\leq n}) + \mu_B((T_{\rm
uns})_{\leq{n}})}{\mu_A((\widehat{A})_{\leq n}) +
\mu_B((T)_{\leq{n}})}\right) \cdot (1 - \varepsilon)^{k-1}
\textrm{ for odd $k$}$$ and
$$\rho_{\mu}^{n,k}(\Sigma_{\rm uns},\Sigma) <
 \left(\fracd{\mu_A((\widehat{A})_{\leq n}) + \mu_A((S_{\rm
uns})_{\leq{n}})}{\mu_A((\widehat{A})_{\leq n}) +
\mu_A((S)_{\leq{n}})}\right) \cdot (1 - \varepsilon)^{k-1}
\textrm{ for even $k$}.$$

In both cases the limit of frequencies
$\rho_{\mu}^{n,k}(\Sigma_{\rm uns},\Sigma)$ while $d(n,k)
\rightarrow \infty$ exists, equal to zero and doesn't depend on a
particular choice of this direction. Moreover, for large enough
$n$ and $k$ these frequencies bounded by $(1-\varepsilon)^{k-1}$
and this completes the prove.

%
%

The proof of the theorem for $\mathcal{CRF}$ is
straightforward.\eproof

\subsection{Regular and stable normal forms: $L-$measure and
$\lambda_L-$measure}\label{subsection:lamda_l}

Now we are ready to formulate and prove one of the most important
results of this work.

\begin{maintheorem}\label{nf_lambdaLmeasurable} Let $G=\mathop{A\ast B}\limits_C$ be an amalgamated product, where
$A,B,C$ are free groups of finite rank. If $C$ of infinite index
either in $A$ or in $B,$ then sets of all unstable
$\mathcal{NF}_{\rm uns}$ and all singular $\mathcal{NF}_{\rm sin}$
normal forms are exponentially
$\lambda_{\mathcal{NF}}-$measurable, where $\mathcal{NF} = \{
\mathcal{EF}, \mathcal{RF}, \mathcal{CNF}, \mathcal{CRF} \}$.
\end{maintheorem}
\proof Observe, that as in Theorem \ref{nf:basic} it is sufficient
to show the result only for unstable normal forms. Suppose $C$ has
an infinite index in $A$.

Consider the set of all freely-reduced normal forms first. Since
$C$ has an infinite index in $A$, there exists at least one stable
representative $s \in S_{\rm st}$. Set $W = \{ y_i s y_j | y_i,
y_j \in Y \cup Y^{-1} \}$; and let, as in Lemma
\ref{le:fp_Lcones}, notation $F_0$ mean the set of all words in
$F=A\ast B$, that doesn't contain any element of $W$ as a subword,
and let $F_1$ be the supplement of $F_0$ in $F$.

Due to Lemma \ref{nf_structure}, the set of all unstable
freely-reduced forms can be written as $\mathcal{EF}_{\rm uns} =
(\mathop{\cup}\limits_{s \in S_{\rm uns}}Cs) \ast
(\mathop{\cup}\limits_{t \in T_{\rm uns}}Ct)$ and therefore
doesn't contain words having subwords of the type $y_i s y_j$,
i.e. $\mathcal{EF}_{\rm uns} \subseteq F_0$. The set $\mathcal{EF}
= A \ast B$ is a free group and all cones $C(u)$ in $A \ast B$ are
precisely all reduced words in this group that start from $u$.
Obviously, all such cones $C(u)$  have nontrivial intersection
with $F_1$, and by Lemmata \ref{nf_structure} and
\ref{le:fp_Lcones} it follows that unstable freely-reduced forms
are exponentially $\lambda_{\mathcal{EF}}-$measurable.

Let us consider the set of all nontrivial reduced forms and the
subset of all unstable forms in it. By Lemma \ref{nf_structure} we
have: $\mathcal{RF}_{\rm uns} = (A_1' \setminus C) \ast (B_1'
\setminus C)$, where $A_1' = \mathop{\cup}\limits_{s \in S_{\rm
uns}}Cs \,\,\textrm{ and } \,\,B_1' = \mathop{\cup}\limits_{t \in
T_{\rm uns}}Ct$. Therefore, the set of all unstable forms is
regular and $\mathcal{RF}_{\rm uns} \subseteq F_0$.

The set of all reduced forms $\mathcal{RF} = (A \setminus C) \ast
(B \setminus C)$ is regular and prefix-closed in $A\ast B$ by
lemma \ref{nf_structure}, and since $C \cap F_1 = \emptyset$, the
set $\mathcal{RF}$ has nonempty intersection with $F_1$. Then by
Lemma \ref{le:fp_Lcones} the set $\mathcal{RF}_{\rm uns}$ is
exponentially $\lambda_{\mathcal{RF}}-$measurable.

Let us prove the theorem for $\mathcal{CNF}$ given by formulae
(\ref{nf_cnf2}) and (\ref{nf_cnf_uns}):

$\mathcal{CNF} = (A \circ (S \ast T)_1) \sqcup (B \circ (S \ast
T)_2)$ and  $\mathcal{CNF}_{\rm uns} = (A \circ (S_{\rm uns} \ast
T_{\rm uns})_1) \sqcup (B \circ (S_{\rm uns} \ast T_{\rm
uns})_2).$

As in Theorem \ref{nf:basic} above, we prove the statement of this
theorem for the pair $\Sigma=(A \circ (S \ast T)_1)$ and
$\Sigma_{\rm uns}=(A \circ (S_{\rm uns} \ast T_{\rm uns})_1)$.

The set $\Sigma_{\rm uns}$ is a (regular) subset of $F_0$. To show
that $\Sigma_{\rm uns}$ doesn't contain non-small $\Sigma-$cones,
decompose both sets using (\ref{free_product}). Then

$\Sigma_{\rm uns} = A (T_{\rm uns} \cup T_{\rm uns}S_{\rm
uns}T_{\rm uns} \cup T_{\rm uns} (S_{\rm uns}T_{\rm uns})^{\ast}
S_{\rm uns})$ and $\Sigma = A (T \cup TST \cup T (ST)^{\ast} S)$.

All cones in $\Sigma$, except, may be, $\Sigma-$cones from $AT$,
have nonempty intersection with $F_1$ and so they can't be
contained in $\Sigma_{\rm uns}$; but $AT$ doesn't have non-small
$\Sigma-$cones, and therefore, the set $\mathcal{CNF}_{\rm uns}$
is exponentially $\lambda_{\mathcal{CNF}}-$measurable.

The proof of the theorem for $\mathcal{CRF}$ is
straightforward.\eproof


\subsection*{Acknowledgements} The authors thank Alexandre Borovik for inspiration and useful
discussions.

\bigskip

\normalsize

\vfill

\noindent \textsf{Elizaveta Frenkel, Moscow State University,
GSP-1, Leninskie gory, 119991, Moscow, Russia}

\noindent {\tt lizzy.frenkel@gmail.com}

\medskip

\noindent \textsf{Alexei G. Myasnikov, Schaefer School of
Engineering and Science, Department of Mathematical Sciences,
Stevens Institute of Technology, Castle Point on Hudson, Hoboken
NJ 07030-5991, USA. }

\noindent {\tt amiasnikov@gmail.com}

\medskip
\noindent \textsf{Vladimir N. Remeslennikov, Omsk Branch of
Mathematical Institute SB RAS, 13 Pevtsova Street, Omsk 644099,
Russia}

\noindent {\tt remesl@iitam.omsk.net.ru}

\end{document}